\documentclass[12pt,english]{amsart}
\usepackage{amsmath,amssymb}
\usepackage{bm,color,hyperref}
\usepackage{graphicx} 
\usepackage{tikz}
    \usetikzlibrary{cd}
\usepackage{standalone}
\usepackage{enumitem} 

\usepackage{tikz}
 \usetikzlibrary{positioning}
 \usetikzlibrary{decorations.markings}
 \usetikzlibrary{arrows.meta}
\usetikzlibrary{decorations.pathreplacing}
\usepackage{longtable,multirow}

\theoremstyle{plain}
\newtheorem{thm}{Theorem}[section]
\newtheorem{lem}[thm]{Lemma}

\newtheorem{prop}[thm]{Proposition}
\newtheorem{conj}{Conjecture}
\newtheorem{quest}{Question}

\theoremstyle{definition}
\newtheorem{defn}[thm]{Definition}

\theoremstyle{remark}
\newtheorem{ex}{Example}[section]

\let\emptyset\varnothing    

\title{Sierpi\'{n}ski Gaskets Julia set and Postcritically finite fractals}
\begin{author}[S. Silvestri]{Stefano Silvestri}
\email{silvestri.dev@gmail.com}
\address{Rome, Italy}
\end{author}

\begin{document}
\begin{abstract}
Devaney's work on the Julia sets of the McMullen family, established these rational maps as central objects in holomorphic dynamics. For Misiurewicz parameters, that is when the orbits of the critical points are all strictly preperiodic to a repelling cycle, the topological point of view taken by Kameyama links the dynamics of such rational maps to that on the attractors of invertible self-similar iterated function system. Ero\u{g}lu, Rohde, and Solomyak showed that in favorable cases the topological conjugacy is in fact quasisymmetric. We survey the machinery underlying these results and propose directions for further research. 
\end{abstract}

\maketitle

\section{Introduction}
Much of Robert Devaney's later work was focused on the study of the numerous topological objects arising from a family of holomorphic functions called the McMullen maps. This family, first studied by McMullen \cite{McMullen1988}, is defined as follows: for a fixed pair $(n,m)$ of integers with $n\geq2$ and $m\geq1$, define for $c\in\mathbb{C}$ the rational map \[R_{n,m,c}(z)=z^n+\frac{c}{z^m}.\]

The topological objects studied by Devaney appear both in the dynamical and the parameter plane (the $c$-plane) and include, among many others,
\begin{itemize}
    \item[-] Sierpi\'{n}ski curves, which are planar sets that are homeomorphic to the well-known Sierpi\'{n}ski carpet fractal;
    \item[-] generalized Sierpi\'{n}ski gaskets, i.e. sets obtained by a similar recursive process defined for the Sierpi\'{n}ski triangle but applied to the closed unit disk, as starting set, and by removing polygons;
    \item[-] Cantor necklaces, sets homeomorphic to the Cantor middle-thirds set with the removed open intervals replaced by open disks.
\end{itemize} 

The rational maps $R_{n,m,c}$ are special in the sense that each possesses a single critical orbit up to symmetry. Excluding the superattracting fixed point at $\infty$ and the pole at $0$, each map $R_{n,m,c}$, for fixed $n$ and $m$, possess $n+m$ critical points that behave symmetrically with respect to the angle $2\pi/(n+m)$. 

The focus of this article is on those parameters $c$ for which the orbits of the critical points remain bounded and are all strictly preperiodic to a repelling cycle. For these parameters the Julia set is a generalized Sierpi\'{n}ski gasket. The goal is to highlight a connection between these maps and the linear setting of iterated function systems of contracting similarities of the plane.

\begin{figure}
    \centering
    \includegraphics[width=0.4\linewidth]{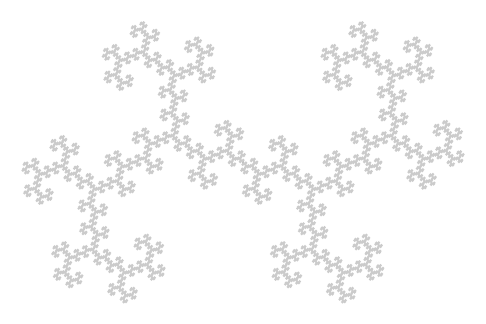}
    \includegraphics[width=0.4\linewidth]{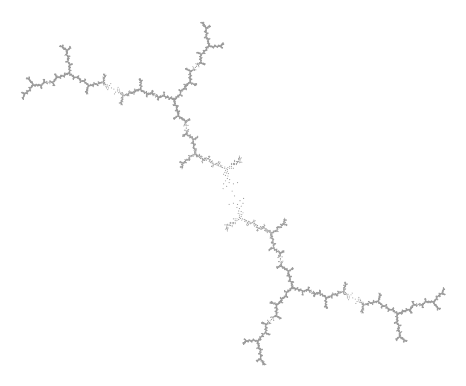}

    \includegraphics[width=0.4\linewidth]{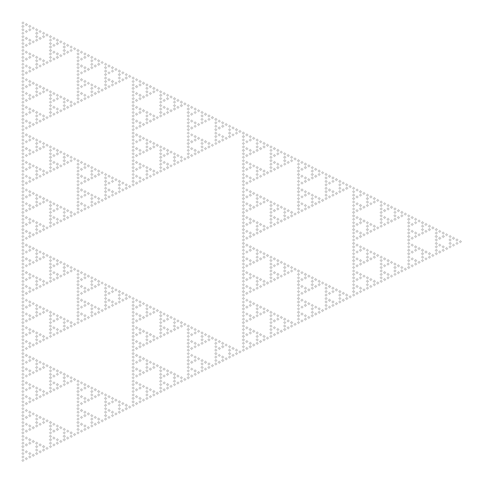}
    \includegraphics[width=0.4\linewidth]{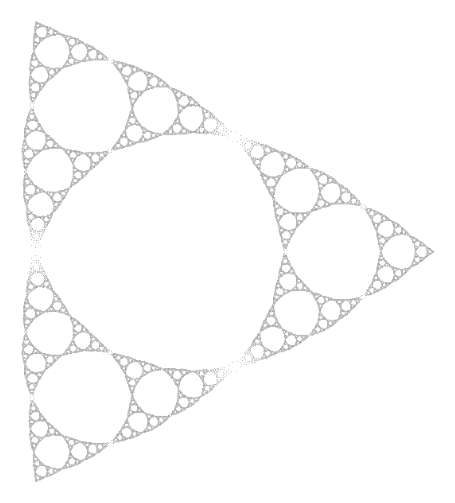}
    
    \includegraphics[width=0.4\linewidth]{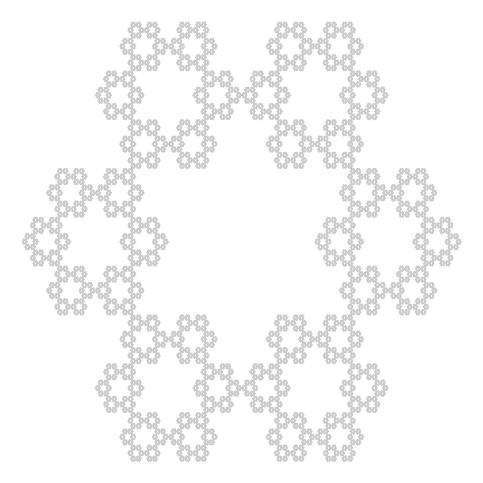}
    \includegraphics[width=0.4\linewidth]{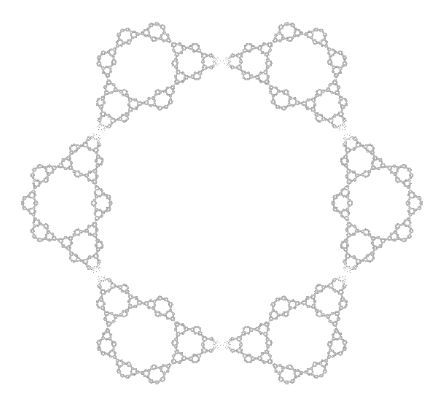}
    \caption{Examples of known attractors of IFS with $2,3$, and $6$ maps which have a corresponding Julia set in the quadratic or the McMullen family.}
    \label{fig:attractorvsJulia}
\end{figure}

Many familiar self-similar sets, such as the Sierpi\'{n}ski triangle, arise as attractors of finite iterated function systems (IFS) consisting of maps of the form $f_{k}(z)=\lambda_k z+v_k$ with $\lambda_k\in\mathbb{D}$ and $v_k\in\mathbb{C}$. The attractor $A$ of an IFS $\{f_k\}_{k=0}^N$ is the unique non-empty compact set satisfying \[A=\bigcup_{k=0}^{N-1}f_k(A).\] When the pieces $f_k(A)$ overlap in a sufficiently controlled way, the IFS can be ``inverted'': suppose that for $i\neq j$ the sets $f_i(A)\cap f_j(A)$ are finite and $f_i^{-1}(x)=f_j^{-1}(x)$ for each $x\in f_i(A)\cap f_j(A)$, then there exists an expanding map $q$ on $A$ such that $q=f_j^{-1}$ when restricted to $f_j(A)$. The resulting dynamical system $(A,q)$ is topologically conjugate to a rational map acting on its Julia set. For example, 
\begin{itemize}
    \item[-] when $N=2$ and $f_0(A)\cap f_1(A)=\{0\}$ then there are values of $\lambda_0$ and $\lambda_1$ for which the system $(A,q)$ is conjugate \cite{Bandt2002} to $(J_c,z^2+c)$ for some Misiurewicz parameter $c$;
    \item[-] when $N=3$ and $A$ is the Sierpi\'{n}ski triangle, then $(A,q)$ is conjugate \cite{Kameyama2000} to the system $(J_c,R_{2,1,c})$ for $c=-16/27$;
    \item[-] when $N=6$ and $A$ is the Hexagasket, $(A,q)$ is conjugate \cite{ErogluRohdeSolomyak2010} to $(J_c,R_{4,2,c})$ for $c=4/9\exp(i\pi/3)$.
\end{itemize}
A plot placing side by side the attractor and the Julia set is shown in Figure \ref{fig:attractorvsJulia}.
Note that each of the attractor in the examples above can be obtained by various IFS, e.g. the Sierpi\'{n}ski triangle is the attractor of the IFS $\{0.5 z+\xi_3^k\}_{k=0}^2$ and of $\{0.5\xi_3^{2k} z+\xi_3^k\}_{k=0}^2$ where $\xi_3=\exp(i2\pi/3)$, but the IFS is ``invertible'' only in the latter case. We revisit these example with more details in Section \ref{subsubsect:knownexamples}. 

Ero\u{g}lu, Rohde, and Solomyak \cite{ErogluRohdeSolomyak2010} showed that the topological conjugacy is quasisymmetric under the hypotheses that $A$ has bounded turning and the target rational map has no parabolic cycles. A connected and locally connected set $X\subset \mathbb{C}$ is said to have bounded turning if for every pair of points $x,y\in X$ there exists a continuum $K\subset X$ containing both points with diameter at most a fixed multiple of $\lvert x-y\rvert$. In particular, bounded turning sets do not have cusps but can exhibit fractal behavior. A homeomorphism $f$ between two metric spaces $(X,d_X)$ and $(Y,d_Y)$ is quasisymmetric if there is a homeomorphism $\eta:[0,\infty)\to[0,\infty)$ such that for every $t>0$ and $a,b,x\in X$ \[d_X(x,a)\leq t d_X(x,b)\implies d_Y(f(x),f(a))\leq \eta(t)d_Y(f(x),f(b)).\]
A quasisymmetric conjugacy respects the metric geometry of the two spaces rather than only their combinatorial structure. Ha\"{i}ssinsky and Pilgrim, within the theory of coarse expanding conformal (CXC) dynamical systems, obtained the quasisymmetric conjugacy as a consequence of a general uniformization theorem for metrically CXC systems. 

Related is also the work of Qiu, Yang, and Yin \cite{QiuYangYin2018} who studied the quasiconformal geometry and gave a quasisymmetric classification of the Julia sets of McMullen maps for $n,m\geq2$ satisfying $1/n+1/m<1$ when the Julia set is a Sierpi\'{n}ski carpet. The dynamics on Sierpi\'{n}ski sets has also recently attracted the attention the field of people in the field of geometric group theory \cite{LodgeLyubichMerenkovMukherjee2023,BelkForrest2025} and quasiconformal geometry \cite{Ntalampekos2019,LuoZhang2025,LuoNtalampekos2026}.

The correspondence between IFS attractors and Julia sets of rational maps has so far been established example by example, without a systematic account of which self-similar sets correspond to which Julia sets. The tools required for such an account are available in the literature but have not been previously assembled.

This survey has two aims. First, we review the relevant tools, organize the known instances of the correspondence, and make explicit the mechanism that produces the quasisymmetric conjugacy. Second, we propose to extend this mechanism beyond the low-degree examples currently documented, to a dictionary linking the ``fractal $N$-gons'' to the McMullen family of degree $n+m=N$. We describe what such dictionary should look like, the evidence supporting it, and the obstacles proving it, including connections to recent objects introduced by in \cite{TetenovYudinDrozdov2026} and \cite{CalegariWalker2026}. In the interest of a self-contained account, we do not delve into the CXC dynamics.

\section{Iterated Function Systems}
\label{sect:IFS}
Let $\lambda\in \mathbb{D}$ and fix an integer $N\geq2$. For $0\leq j<N$ define the IFS $\Phi_N:=\{f_j\}_j$ of contraction similarities of $\mathbb{C}$:
\[
    f_j(z)=\xi_N^j+\lambda z,\quad \xi_N=\exp(2\pi i/N).
\]
Standard arguments show that there is a unique invariant non-empty compact set $A_\lambda\subset\mathbb{C}$ associated to $\Phi_N$ such that the following equation is satisfied \[A_\lambda=f_0(A_\lambda)\cup f_1(A_\lambda)\cup\cdots\cup f_{N-1}(A_\lambda).\] In other words, $A_\lambda$ is the unique fixed point of the function $\Phi_N(X):=\bigcup_{0\leq j<N}f_j(X)$ on the space of compact sets of $\mathbb{C}$. For the rest of the article we will use $\Phi_N(A_\lambda)$ and $A_\lambda$ interchangeably, using the former if we want to stress the number of maps or distinguishing between IFSs. 

The attractor then consists of similar copies $A_k=f_k(A_\lambda)$ of itself and each piece $A_k$ is made up of even smaller copies $A_{kj}=f_k(f_j(A_\lambda))$. It follows that for any integer $n\geq0$ we can consider the set $\Sigma_N^n$ of words $\mathbf{w}=w_0\ldots w_{n-1}$ from the alphabet $\Sigma_N:=\{0,1,\cdots,N-1\}$. Writing $f_{\mathbf{w}}=f_{w_0}\circ\cdots\circ f_{w_{n-1}}$ and $A_{\mathbf{w}}=f_{\mathbf{w}}(A_\lambda)$, we have \[A_\lambda=\bigcup_{n\geq0}\{A_{\mathbf{w}}~\mid~\mathbf{w}\in\Sigma_N^n\}.\]
If we let $n$ tend to infinity, this induces a continuous map $\pi_\lambda:\Sigma_N^\infty\to A_\lambda$, called the \emph{address map}, from the set of infinite (right) sequences $\Sigma_N^\infty$ onto the attractor.
Another characterization of the attractor is then given by \[A_\lambda=\pi_\lambda(\Sigma^\infty_N)=\bigg\{\sum_{j=0}^\infty\xi_N^{ \varepsilon_j}\lambda^j~\bigm\vert~\varepsilon_j\in\Sigma_N\bigg\}.\] 


\subsection{Properties of the attractor}
\label{subsect:propertiesAttractor}
The simple definition of $\Phi_N(A_\lambda)$ implies that the sets $f_j(A_\lambda)$ are translations of each other: for any $k,j\in\Sigma_N$ we have \[A_k=\lambda A_\lambda +\xi_N^k=\lambda A_\lambda +\xi_N^k-\xi_N^j+\xi_N^j=A_j+\xi_N^k-\xi_N^j.\] Repeating the above argument for pieces at level $n$ we obtain 
\begin{lem}
\label{lem:piecestranslation}
    For any $n\geq0$ and for each distinct finite word $\mathbf{w},\mathbf{v}\in\Sigma_N^n$ we have $A_{\mathbf{v}}=f_{\mathbf{v}}f_{\mathbf{w}}^{-1}(A_{\mathbf{w}})$.
\end{lem}

Note that any similarity map, i.e. a change of Cartesian coordinate system, does not change the structure of $A_\lambda$
\begin{lem}
\label{lem:geometricalsimilarity}
    Let $A_\lambda$ be the attractor of the IFS $\Phi_N$. For any similarity map $h:\mathbb{C}\to\mathbb{C}$ the set $B=h(A_\lambda)$ is the self-similar set generated by the IFS $\{g_k=h\circ f_k\circ h^{-1}\}_{k=0}^{N-1}$.
\end{lem}
\begin{proof}
    \[h(A_\lambda)=h\bigg(\bigcup_{k=0}^{N-1}f_k(h^{-1}(B))\bigg)=\bigcup_{k=0}^{N-1}(h\circ f_k\circ h^{-1})(B)=\bigcup_{k=0}^{N-1}g_k(B)=B,\] where the last equality is due to the uniqueness of the attractor of an IFS.
\end{proof}
Therefore, using the IFS $\Phi_N$ is as good as using $\{g_k\}_{k=0}^{N-1}$ since the attractor $A_\lambda$ and $B$ are geometrically similar copies. Observe also that for certain non-trivial choices of $h$, we can actually have $B=A_\lambda$.
\begin{lem}
    \label{lem:rotationalsymmetry}
    Let $s_j(z)=\xi_N^jz$, then $s_j(A_\lambda)=A_\lambda$ where $A_\lambda=\Phi_N(A_\lambda)$. Furthermore, $A_\lambda$ is the attractor of the IFS $\Phi_{N,\sigma}:=\{f_k\circ s_{\sigma(k)}\}_{k=0}^{N-1}$ for any function $\sigma$ on $\Sigma_N$.
\end{lem}
\begin{proof}
    For a fixed $j\in\Sigma_N$, direct computation shows $f_k\circ s_j=s_j\circ f_{k-j}$ for every $k\in\Sigma_N$, proving the first claim.
    
    Let $\sigma:\Sigma_N\to\Sigma_N$ be any function. Since $s_{\sigma(k)}(A_\lambda)=A_\lambda$ for any $k\in\Sigma_N$, then 
    \[A_\lambda=\bigcup_{k=0}^{N-1}f_k(A_\lambda)=\bigcup_{k=0}^{N-1}f_k(s_{\sigma(k)}(A_\lambda))=\bigcup_{k=0}^{N-1}(f_k\circ s_{\sigma(k)})(A_\lambda).\]
\end{proof}
The first statement in the above lemma gives us the $N$-rotational symmetry of the attractor. The second claim tells us that we can modify the defining maps of the IFS without changing the attractor. This is particularly useful in determining whether we can define an expanding map on the attractor, see Section \ref{sect:IIFSAndConjugacy}.

There is a generally accepted condition that controls the overlaps $A_j\cap A_k$ of an attractor due to Mor\'{a}n \cite{Moran1946}. The attractor is said to satisfy the \emph{open set condition}, abbreviated OSC, if there exists nonempty open set $V\subset\mathbb{C}$ with \[\bigcup_{k=0}^{N-1}f_k(V)\subset V\text{ and }f_j(V)\cap f_k(V)=\emptyset\text{ for }j\neq k.\]
For a given IFS, however, finding such open set $V$ is not always an easy task. Equivalent definitions have been found in order to get around this issue. 

\begin{thm}
\label{thm:OSC}
Let $A$ be attractor of an IFS $\{f_j\}_{j=1}^m$ with contraction factors $r_j$. The following are equivalent
 \begin{enumerate}
    \item The IFS satisfies the OSC.
    \item \underline{Positivity of the $\alpha$-Hausdorff measure} {\rm (Mor\'{a}n \cite{Moran1946} and Schief \cite{Schief1994})}: $\mu^\alpha(A)>0$ where $\alpha$ is the similarity dimension of $A$ given by $\sum_j r_j^\alpha=1$. 
    \item \underline{The finite clustering property} {\rm (Schief \cite{Schief1994})}: there exists an integer $N$ such that for every $A_{\mathbf{w}}$ of diameter $\varepsilon$, for a finite length word $\mathbf{w}$ of symbols in $\Sigma_m$, there are at most $N$ incomparable pieces $A_{\mathbf{v}}$ of diameter $\geq\varepsilon$ that can intersect the $\varepsilon$-neighborhood of $A_{\mathbf{w}}$. We say that $A_{\mathbf{w}}$ and $A_{\mathbf{v}}$ are incomparable if $\mathbf{w}$ is not a prefix $\mathbf{v}$ and viceversa.
    \item \underline{The neighbor map condition} {\rm (Bandt--Graf \cite{BandtGraf1992})}: the identity map is not an accumulation point of the set of neighbor maps of $A$. A neighbor map has the form $h=f_{\mathbf{w}}^{-1}f_{\mathbf{v}}$ where $\mathbf{w}=w_0\cdots w_k$ and $\mathbf{v}=v_0\cdots v_k$ with $w_0\neq v_0$ for some finite $k>0$. 
 \end{enumerate}
\end{thm}

A significant consequence of OSC for our purpose is the following theorem. If we define the equivalence relation $\sim$ on $\Sigma_N^\infty$ as $\mathbf{w}\sim\mathbf{v}$ if $\pi(\mathbf{w})=\pi(\mathbf{v})$, then 
\begin{thm}[\cite{BandtMesing2009}]
    If the open set condition holds, all equivalence classes of addresses are finite.
\end{thm}
Thankfully in our setting OSC is satisfied since
\begin{thm}[\cite{BandtRao2007}]
    \label{thm:FiniteOverlapImpliesOSC}
    Let $A$ be a connected self-similar set in the plane. If $A_j\cap A_k$ is a finite set for $i\neq j$, then OSC holds.
\end{thm}


\subsection{Connectedness locus of the attractor and its subsets}
\label{subsect:mandelbrotsetsforIFS}
Connectedness of the attractor is a well-studied property of IFS, especially since, for self-similar sets, it implies local connectedness and arcwise connectedness \cite{BandtKeller1991,Hata1985}. A general theorem of Hata \cite{Hata1985} states that for self-similar sets, either the attractor is connected or it is homeomorphic to a Cantor set. In our case, by the rotational symmetry if $A_0\cap A_1$ is non-empty then so are $A_k\cap A_{k+1}$ for any $0\leq k<N$. Hence, looking at the possible intersections of the pieces $A_0$ and $A_1$ is enough.
\begin{thm}[\cite{BandtHung2008b}]
\label{thm:connectednesscriterion}
    The attractor $\Phi_N(A_\lambda)$ is connected if and only if $A_0\cap A_1\neq\emptyset$. If $\Phi_N(A_\lambda)$ is disconnected then it is totally disconnected, and all pieces $A_k$ are disjoint.
\end{thm}
Define the connectedness locus of the IFS 
\[\mathbb{B}_N=\{\lambda\in\mathbb{D}~\mid~\Phi_N(A_\lambda)\text{ is connected}\}.\]
Using the power series representation of a point $z\in A_\lambda$ and Theorem \ref{thm:connectednesscriterion} we obtain another characterization of the above set as
\[\mathbb{B}_N=\{\lambda\in\mathbb{D}~\mid~\sum_{j\geq0}\varepsilon_j\lambda^j=0,~\varepsilon_j\in \Delta_N\}\] where $\Delta_N:=\{(\xi_N^a-\xi_N^b)/(1-\xi_N)~\mid~a,b\in\Sigma_N\}$. 
This set is a natural object to study and it can be thought of as the linear analogue of the connectedness locus for Julia sets of rational maps.

Barnsley and Harrington \cite{BarnsleyHarrington1985} in 1985 started the investigation of $\mathbb{B}_2$, and, in fact, introduced the idea of a connectedness locus for a family of IFS. They proved basic geometric properties of it and started the investigation of $\partial \mathbb{B}_2$. Between 1988 and 1992 Bousch \cite{Bousch1988,Bousch1993} proved $\mathbb{B}_2$ is connected and locally connected. In 2002 Bandt \cite{Bandt2002} proved the existence of a hole in $\mathbb{B}_2$, aside from the central one, and conjectured that there are in fact infinitely many holes. In 2017 Calegari, Koch, and Walker \cite{CalegariKochWalker2017} gave a positive answer to the conjecture by showing the existence of infinite sequences of holes accumulating at parameters in $\partial \mathbb{B}_2$. These are buried points of $\partial \mathbb{B}_2$, that is not accessible from $\mathbb{D}\setminus \mathbb{B}_2$. In \cite{CalegariKochWalker2017} the authors also settled another conjecture by Bandt that the interior points of $\mathbb{B}_2$ are dense away from the real axis. In 2005 Solomyak \cite{Solomyak2005} proved that for points on $\partial \mathbb{B}_2$, there is an asymptotic similarity between $\mathbb{B}_2$ and the attractor of a related IFS$\{-1+\lambda z, \lambda z,1+\lambda z\}$. In 2023, P\'erez and Silvestri \cite{SilvestriPerez2023} exploiting the result of Solomyak, gave sufficient conditions for $\lambda\in \mathbb{B}_2$ to be path accessible from $\mathbb{D}\setminus \mathbb{B}_2$. In the latest article by Calegari and Walker \cite{CalegariWalker2026} it is shown that aside from the points $\pm i/\sqrt{2}$, the entire annulus of minimum radius $1/\sqrt{2}$ inside the unit disk is contained in the interior of $\mathbb{B}_2$. Consequently, $\partial\mathbb{B}_2$ is contained in the annulus with radii $1/2$ and $1/\sqrt{2}$.

The study of $\mathbb{B}_N$ for $N>2$ has also been initiated. Bandt and Hung \cite{BandtHung2008b} gave elementary properties of $\mathbb{B}_N$ and showed $\mathbb{B}_N$ is regular closed, i.e. that $\mathbb{B}_N=\overline{\mathbb{B}_N^\circ}$, for $N\geq3$ but not $N=4$. Himeki and Ishii \cite{HimekiIshii2020} filled the gap by extending the proof of \cite{CalegariKochWalker2017}. Nakajima \cite{Nakajima2024} extended Bousch's arguments and gave a broad criterion for the set of roots of power series with prescribed coefficients to be connected and locally connected. Some rigidity result on the set of roots of power series with coefficients from an arbitrary finite set are given in \cite{Kewart2025}. In \cite{BandtHung2008b} the authors also provided sharp bounds for $\mathbb{B}_N$ and simplified the known formula for the minimum value of $\lambda$ that produces ``just touching'' attractors when $N\not\equiv0\bmod4$.

Let us return back to the analogy with the parameter plane for rational map. When $N=2$, it seems natural to make the parallel to the classical Mandelbrot set \[\mathcal{M}=\{c\in\mathbb{C}~\mid~J_c=J(z^2+c)\text{ is connected}\}.\] 
As we discussed above, there are properties of $\mathcal{M}$ that have analogies to $\mathbb{B}_2$: connectedness, local connectedness (even though it is just conjectured for $\mathcal{M}$), asymptotic similarity at points on $\partial\mathcal{M}$. However, the correspondence breaks down since $\mathbb{B}_2$ has holes while the Mandelbrot set does not. There is, in fact, another contrasting property between the two parameters planes: the Mandelbrot set has an alternative definition \[\mathcal{M}=\{c\in\mathbb{C}~\mid~0\text{ is in the filled Julia set}\}\] which has a direct counterpart in the IFS setting to \[\mathbb{B}_{2,0}=\{\lambda\in\mathbb{D}~\mid~0\in \Phi_2(A_\lambda)\}\] but, alas, $\mathbb{B}_{2,0}$ is strictly a subset of $\mathbb{B}_2$. The fastest way to see this inclusion is to rely on their definition in terms of roots of power series. For a given set $S\subset\mathbb{C}$ define the set of power series with coefficient in $S$ as $\mathcal{P}_{S}:=\{F(z)=1+\sum_{k\geq1}\varepsilon_jz^j~\mid~\varepsilon_j\in S\}$. It follows by Theorem \ref{thm:connectednesscriterion} \[\mathbb{B}_2=\big\{\lambda\in\mathbb{D}~\mid~\exists F\in \mathcal{P}_{\{-1,0,1\}}\text{ such that }F(\lambda)=0\big\};\]
but having $0\in A_0\cap A_1$ is a stronger restriction and implies that \[\mathbb{B}_{2,0}=\big\{\lambda\in\mathbb{D}~\mid~\exists F\in \mathcal{P}_{\{-1,1\}}\text{ such that }F(\lambda)=0\big\}.\]
By Lemma \ref{lem:rotationalsymmetry} the attractor $\Phi_2(A_\lambda)$ is symmetric about the origin, and $0$ is also the center of translational symmetry for $\Phi_2$. For these reasons, $0$ is a sort of critical point for the attractor, and consequently the reason that $\mathbb{B}_{2,0}$ is relevant.

Much of the properties of $\mathbb{B}_2$, such as connectedness \cite{Bousch1993} and existence of holes \cite{CalegariKochWalker2017}, extend to $\mathbb{B}_{2,0}$, but the more striking is the result of Tiozzo \cite{Tiozzo2020} which connects our linear setting to quadratic dynamics: $\mathbb{B}_{2,0}$ equals the closure in $\mathbb{D}$ of the set containing all Galois conjugates of growth rates of unimodal maps which are critically periodic. In other words, $\mathbb{B}_{2,0}$ is the set, restricted to the unit disk, introduced by Thurston \cite{Thurston2014} and studied also by Bray, Davis, Lindsey, and Wu \cite{BrayDavisLindseyWu2021,LindseyWu2023}. Together with Tiozzo, Lindsey and Wu \cite{LindseyTiozzoWu2025} also defined analogues of the Thurston' set for each principal vein in the Mandelbrot set. 

A further subset of $\mathbb{B}_{2,0}$ deserve attention. Define \[\mathcal{T}_2=\{\lambda\in\mathbb{D}~\mid~A_0\cap A_1\text{ is a singleton}\}.\]
Since $\Phi_2(A_\lambda)$ is invariant under $z\mapsto-z$ as a set, a singleton overlap of the two pieces $A_0$ and $A_1$ implies that such point is the origin. When $\lambda\in\mathcal{T}_2$, the attractor $\Phi_2(A_\lambda)$ is a dendrite and $\mathcal{T}_2$ is precisely the parameter set considered in establishing the quasisymmetric conjugacy between $(\Phi_2(A_\lambda),q_\lambda)$ and $(J_c,z^2+c)$ at Misiurewicz parameters $c$.
The relationship between $\mathcal{T}_2$ and the sets discussed above is not yet fully understood. However, Solomyak was able to give an important characterization of this set. 
\begin{lem}[\cite{Solomyak2005}]
    \label{lem:no0fordendrite}
    Let $S=\{-1,0,1\}$ then $A_0\cap A_1=\{0\}$ if and only if there exists a unique power series $F\in \mathcal{P}_S$ vanishing at $\lambda$ with no $0$ coefficient.
\end{lem}
The following questions are still open
\begin{quest}
    Is $\mathcal{T}_2=\mathbb{B}_{2,0}\cap \partial \mathbb{B}_2$?
\end{quest}
\begin{quest}
    Is $\mathcal{T}_2$ a Cantor set in $\mathbb{B}_2$?
\end{quest}
Using Rouche's Theorem and careful estimates Solomyak \cite[Theorem 2.4]{Solomyak2005} was able to prove that there exist uncountably many $\lambda\in\mathcal{T}_2$.
\begin{quest}
    Is $\mathcal{T}_2$ dense in $\partial\mathbb{B}_2$?
\end{quest}
\begin{quest}
    Is every $\lambda\in\mathcal{T}_2$ accessible from $\mathbb{D}\setminus \mathbb{B}_{2,0}$ only through the central hole of $\mathbb{B}_{2,0}$? 
\end{quest}
The reason to use $\mathbb{B}_{2,0}$ rather than $\mathbb{B}_2$, in the last question is due to the example of a buried parameter found in \cite{CalegariKochWalker2017}, which is shown to be in $\mathbb{B}_{2,0}\cap\partial\mathbb{B}_2$.

Given the goal of this survey, one may ask: is there a set $\mathbb{B}_{N,0}$ for $N>2$? What about $\mathcal{T}_N$? For the latter, the definition of $\mathcal{T}_2$ in terms of a property of the attractor, as opposed to the one that uses power series, could in principle be reused. The set $\mathbb{B}_{N,0}$ is also delicate: what should play the role of $0$? In either case, characterization using power series seems difficult to find.
As far as we know, these sets have not been defined in the literature. Here we propose a possible definition inspired by the aim of this article. The proved and conjectured properties of $\mathbb{B}_{2,0}$ and $\mathcal{T}_2$ become, then, open questions for these newly defined set $\mathbb{B}_{N,0}$ and $\mathcal{T}_N$.

The origin is the center of symmetry of $\Phi_2(A_\lambda)$, both translationally and rotationally, so it acts as a sort of critical point; for $N>2$ we no longer have such a point. Let $N>2$ and for $0\leq k<N$ define the following half lines in the plane 
\[L_k:=\{t(1+\xi_N)\xi_N^k~\mid~t\geq0\}\]
then \[\mathbb{B}_{N,0}:=\{\lambda\in\mathbb{D}~\mid~(A_0\cap A_1)\cap L_0\neq \emptyset\}\]
and \[\mathcal{T}_N:=\{\lambda\in\mathbb{D}~\mid~A_0\cap A_1=\{p\}\text{ and }p\in L_0\}.\]
When $N=2$, in fact, the lines $L_0$ and $L_1$ are degenerate and coincide with the origin. The direction of $L_k$ is the same as the midpoint between $\xi_N^k$ and $\xi_N^{k+1}$. Consequently, if $\lambda\in\mathbb{R}$ is such that $A_0$ and $A_1$ just touch, the points in the intersection lie on $L_0$. More reasons for these definitions will be given later when discussing the conjugacy between IFS and rational maps.

\begin{figure}
    \centering
    \includegraphics[width=0.45\linewidth]{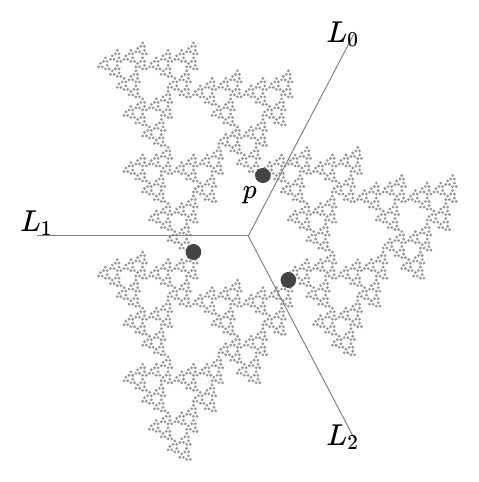}
    \includegraphics[width=0.45\linewidth]{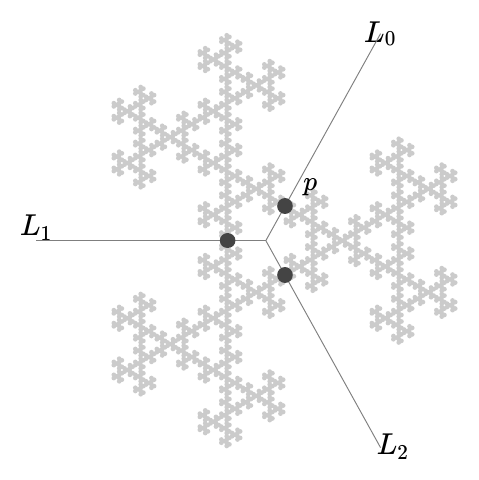}
    \caption{Two examples of just touching IFS when $N=3$. The parameter $\lambda$ for the one on the left does not lie in $\mathcal{T}_3$; while it does for the one on the right. In both cases, there is no $\sigma:\Sigma_3\to\Sigma_3$ for which $(f_0\circ s_{\sigma(0)})^{-1}(p)=(f_1\circ s_{\sigma(1)})^{-1}(p)$.}
    \label{fig:badSierpinskigasket}
\end{figure}

We remark that with our definition of $\mathcal{T}_N$, we are considering a subset of those $\lambda$ for which $A_0\cap A_1$ is a singleton. Consider the example given in \cite[Example 1]{BandtHung2008b} and shown in Figure \ref{fig:badSierpinskigasket}: let $N=3$ and $\lambda=0.5+0.5(2-\sqrt{3})i$ then the attractor satisfies $A_0\cap A_1=\{p\}$ but $p=-\lambda\xi_3^2=\lambda(1+\xi_3)$ does not lie on $L_0$. As it can be deduced from the picture, the addresses of $p$ are $\pi_\lambda^{-1}(p)=\{\mathbf{01}\overline{\mathbf{1}},\mathbf{10}\overline{\mathbf{2}}\}$, where $\overline{x}=xxx\cdots$. Observe that for any function $\sigma$ on $\Sigma_3$ we have that \[(f_0\circ s_{\sigma(0)})^{-1}(p)=(f_1\circ s_{\sigma(1)})^{-1}(p)\iff\frac{\pi_\lambda(\mathbf{1}\overline{\mathbf{1}})}{\pi_\lambda(\mathbf{0}\overline{\mathbf{2}})}=\frac{f_0^{-1}(p)}{f_1^{-1}(p)}=\xi_3^{\sigma(1)-\sigma(0)},\] which one can check is not the case for any value of $\sigma(1)-\sigma(0)$ modulo $3$.
We believe that the reason this IFS fails to be invertible lies in the structure of the two addresses of $p$: if $q$ existed then the two addresses defining $p$ should be ``symmetric'' after the first symbol. More details are in Section \ref{sect:IIFSAndConjugacy}. 
\begin{lem}
\label{lem:ipcfsubsetofT}
    Let $N\geq2$ and suppose $\lambda$ is such that $A_0\cap A_1=\{p\}$. Let $\sigma:\Sigma_N\to\Sigma_N$ be defined as $\sigma(x)=ax+b\bmod N$ for some $a,b\in\Sigma_N$ with $a\neq0$. If the equation $(f_0\circ s_{\sigma(0)})^{-1}(p)=(f_1\circ s_{\sigma(1)})^{-1}(p)$ is satisfied, then $\lambda\in\mathcal{T}_N$.
\end{lem}
\begin{proof}
    Note that the map $(f_k\circ s_{\sigma(k)})^{-1}$ is always an expanding affine transformation. Two \emph{distinct} affine transformations of the plane agree at exactly two points, one of which is the point at $\infty$. It follows that the unique finite point $w$ where $(f_0\circ s_{\sigma(0)})^{-1}$ equals $(f_1\circ s_{\sigma(1)})^{-1}$ has the form 
    \[w=\frac{\xi_N^a-\xi_N}{\xi_N^a-1}\]
    and this point is necessarily on the line $L_0$, i.e. that $w=t(1+\xi_N)$ for some real $t\geq0$. Indeed, if we consider the function $\mu(z)=(z-\xi_N)/(z-1)$, then the unit circle is mapped by $\mu$ to the line $-L_0\cup L_0$. In particular, since for example $\mu(-1)=0.5(1+\xi_N)$, the short arc from $1$ to $\xi$ is sent to $-L_0$ while the long arc to $L_0$, so the claim follows.

    When $N=2$, since $a\neq0$ then $a=1$ and $(f_0\circ s_{b})^{-1}(p)=(f_1\circ s_{1+b})^{-1}(p)$ forces that $p=0$, that is $\lambda\in\mathcal{T}_2$.
\end{proof}
Unfortunately, we can find counterexamples for the other direction: when $N=3$ and $\lambda=-1/2$, the attractor is a relative of the Sierpi\'{n}ski triangle, see Figure \ref{fig:badSierpinskigasket}. $A_0\cap A_1=\{p\}$ with $p=-\xi_3^2/3=(1+\xi_3)/3=\pi_\lambda(\mathbf{0}\overline{\mathbf{01}})=\pi_\lambda(\mathbf{1}\overline{\mathbf{10}})$ but the ratio $\pi_\lambda(\overline{\mathbf{01}})$ over $\pi_\lambda(\overline{\mathbf{10}})$ is not a power of $\xi_3$.


\subsection{The neighbor concept}
\label{subsect:neighborgraphs}
There is an algebraic framework which is applicable on a large class of fractals, including our families. Such framework has been used not only to determine the Hausdorff dimension of the boundary of fractals (see for example \cite{LauNgai2007}), but also to completely describe the topological structure of the fractal in terms of an equivalence relation on $\Sigma_N^\infty$ \cite{BandtMesing2009}.

Let $\mathbf{w},\mathbf{v}\in\Sigma_N^n$ be such that $A_{\mathbf{w}}\cap A_{\mathbf{v}}\neq\emptyset$, then we define a neighbor map as the function $h=f_{\mathbf{w}}^{-1}f_{\mathbf{v}}$. From Lemma \ref{lem:piecestranslation} is easy to see that $h$ is always a translation and $h(A_\lambda)$ is a copy of $A_\lambda$ that intersects it in the same way $A_{\mathbf{w}}$ intersects $A_{\mathbf{v}}$
\[A_{\mathbf{w}}\cap A_{\mathbf{v}}\neq\emptyset\iff f_{\mathbf{w}}^{-1}(A_{\mathbf{w}})\cap f_{\mathbf{w}}^{-1}(A_{\mathbf{v}})\neq\emptyset\iff A_\lambda \cap h(A_\lambda)\neq\emptyset.\]

Neighbor maps are generated recursively: start with the identity map and for each $j,k\in\Sigma_N$ apply repeatedly the interior automorphism $\alpha
_{jk}(h):=f_j^{-1}\circ h\circ f_k$ of the isometry group of $\mathbb{C}$. The identity map describes the trivial type of two equal pieces, and it is always the starting point of the recursion. Note that the new map $\alpha_{jk}(h)=f_{\mathbf{w}j}^{-1}f_{\mathbf{v}k}$ is in fact a \emph{potential} neighbor map since if $h=f_{\mathbf{w}}^{-1}f_{\mathbf{v}}$ is a neighbor map, i.e. $A_{\mathbf{w}}\cap A_{\mathbf{v}}\neq\emptyset$, it is not guaranteed that $A_{\mathbf{w}j}\cap A_{\mathbf{v}k}\neq\emptyset$ for any $j,k\in\Sigma_N$. It turns out (see \cite{BandtMesing2009}) that the condition $\lvert\alpha_{jk}(h)(0)\rvert\leq2/(1-\lvert\lambda\rvert)$ is necessary but not sufficient to determine whether $\alpha_{jk}(h)$ is an actual neighbor map. 

Nevertheless we select all those potential maps $h$ that satisfy $\lvert h(0)\rvert\leq2/(1-\lvert\lambda\rvert)$. We use these maps as the vertex set of a directed graph. The edge between $h$ and $\alpha_{jk}(h)$ is labeled $jk$. We further reduce the graph by repeatedly removing any vertex without any outgoing edge. The resulting graph is called the \emph{neighbor graph} of $\Phi_N$  at $\lambda$ and we denoted it by $NG_\lambda$. We say that $A_\lambda$ is of \emph{finite type} if the graph $NG_\lambda$ is finite. In this case, from each vertex there starts a path which ends in a directed cycle.

\begin{thm}[\cite{BandtMesing2009}]
\label{thm:equivalentaddresses}
    The pairs of equivalent addresses coincide with the label sequences of infinite edge paths in the neighbor graph. Moreover, for a finite type fractal, the equivalence relation $\sim$ of addresses is defined by the regular language on $\Sigma_N^2$ given by the neighbor graph.
\end{thm}
The above theorem explains how the neighbor graph encodes information about intersections between pieces of the attractor: an edge labeled $s_1t_1$ from the root $id$ to a vertex $h$ exists if and only if $A_{s_1}\cap A_{t_1}\neq\emptyset$; thus an infinite path $s_1t_1, s_2t_2,\ldots$ in the graph corresponds to two addresses $\mathbf{s}=s_1s_2\cdots$ and $\mathbf{t}=t_1t_2\cdots$ which describe the same point $\pi_\lambda(\mathbf{s})=\pi_\lambda(\mathbf{t})$ in the attractor.

\begin{figure}
    \centering
    \begin{tikzpicture}
    
    \node [draw, circle] (0) at (0:0) (id) {\small $id$};
    \node [draw, fill=black, circle] (120) at (120:3) (s0a) {};
    \node [draw, fill=black, circle] (60) at (60:3) (s0b) {};
    \node [draw, fill=black, circle] (0) at (0:3) (s0c) {};
    \node [draw, fill=black, circle] (-60) at (-60:3) (s1a) {};
    \node [draw, fill=black, circle] (-120) at (-120:3) (s1b) {};
    \node [draw, fill=black, circle] (-180) at (-180:3) (s1c) {};

    \path[-{Triangle[]}] (id) edge [right] node {$01$} (s0a);
    \path[-{Triangle[]}] (id) edge [left] node {$10$} (s1a);
    
    \path[-{Triangle[]}] (s0a) edge [above] node {$10$} (s0b);
    \path[-{Triangle[]}] (s0b) edge [right] node {$01$} (s0c);
    \path[-{Triangle[]}] (s0c) edge [right] node {$01$} (s1a);
    \path[-{Triangle[]}] (s1a) edge [below] node {$01$} (s1b);
    \path[-{Triangle[]}] (s1b) edge [left] node {$10$} (s1c);
    \path[-{Triangle[]}] (s1c) edge [left] node {$10$} (s0a);
        
\end{tikzpicture}
    \begin{tikzpicture}
    
    \node [draw, circle] (0) at (0:0) (id) {\small $id$};
    \node [draw, fill=black, circle] (-10) at (-10:3) (s0a) {};
    \node [draw, fill=black, circle] (10) at (10:3) (s0b) {};
    \node [draw, fill=black, circle] (100) at (100:3) (s1a) {};
    \node [draw, fill=black, circle] (120) at (120:3) (s1b) {};
    \node [draw, fill=black, circle] (250) at (250:3) (s2a) {};
    \node [draw, fill=black, circle] (230) at (230:3) (s2b) {};

    \path[-{Triangle[]}] (id) edge [bend left=10, above] node {$20$} (s0b);
    \path[-{Triangle[]}] (id) edge [bend right=10, below] node {$10$} (s0a);
    \path[-{Triangle[]}] (id) edge [bend right=10, right] node {$21$} (s1a);
    \path[-{Triangle[]}] (id) edge [bend left=10, left] node {$01$} (s1b);
    \path[-{Triangle[]}] (id) edge [bend left=10, right] node {$12$} (s2a);
    \path[-{Triangle[]}] (id) edge [bend right=10, left] node {$02$} (s2b);
    \draw[-{Triangle[]}] (s0a) to[out=-45, in=45, looseness=10] node[right] {$01$} (s0a);
    \draw[-{Triangle[]}] (s0b) to[out=-45, in=45, looseness=10] node[right] {$02$} (s0b);
    \draw[-{Triangle[]}] (s1b) to[out=70, in=160, looseness=10] node[above] {$10$} (s1b);
    \draw[-{Triangle[]}] (s1a) to[out=70, in=160, looseness=10] node[above] {$12$} (s1a);
    \draw[-{Triangle[]}] (s2b) to[out=190, in=280, looseness=10] node[below] {$20$} (s2b);
    \draw[-{Triangle[]}] (s2a) to[out=190, in=280, looseness=10] node[below] {$21$} (s2a);
    
\end{tikzpicture}
    \begin{tikzpicture}
    
    \node [draw, circle] (0) at (0:0) (id) {\small $id$};
    \node [draw, fill=black, circle] (0) at (0:3) (s0) {};
    \node [draw, fill=black, circle] (60) at (60:3) (s1) {};
    \node [draw, fill=black, circle] (120) at (120:3) (s2) {};
    \node [draw, fill=black, circle] (180) at (180:3) (s3) {};
    \node [draw, fill=black, circle] (240) at (240:3) (s4) {};
    \node [draw, fill=black, circle] (300) at (300:3) (s5) {};

    \path[-{Triangle[]}] (id) edge [bend right=10, below] node {$45$} (s0);
    \path[-{Triangle[]}] (id) edge [bend left=10, above] node {$21$} (s0);
    \path[-{Triangle[]}] (id) edge [bend right=10, right] node {$50$} (s1);
    \path[-{Triangle[]}] (id) edge [bend left=10, left] node {$32$} (s1);
    \path[-{Triangle[]}] (id) edge [bend right=10, right] node {$01$} (s2);
    \path[-{Triangle[]}] (id) edge [bend left=10, left] node {$43$} (s2);
    \path[-{Triangle[]}] (id) edge [bend right=10, above] node {$12$} (s3);
    \path[-{Triangle[]}] (id) edge [bend left=10, below] node {$54$} (s3);
    \path[-{Triangle[]}] (id) edge [bend right=10, left] node {$23$} (s4);
    \path[-{Triangle[]}] (id) edge [bend left=10, right] node {$05$} (s4);
    \path[-{Triangle[]}] (id) edge [bend right=10, left] node {$34$} (s5);
    \path[-{Triangle[]}] (id) edge [bend left=10, right] node {$10$} (s5);
    \draw[-{Triangle[]}] (s0) to[out=-45, in=45, looseness=10] node[right] {$03$} (s0);
    \draw[-{Triangle[]}] (s1) to[out=20, in=110, looseness=10] node[above] {$14$} (s1);
    \draw[-{Triangle[]}] (s2) to[out=80, in=170, looseness=10] node[above] {$25$} (s2);
    \draw[-{Triangle[]}] (s3) to[out=135, in=225, looseness=10] node[left] {$30$} (s3);
    \draw[-{Triangle[]}] (s4) to[out=190, in=280, looseness=10] node[below] {$41$} (s4);
    \draw[-{Triangle[]}] (s5) to[out=250, in=340, looseness=10] node[below] {$52$} (s5);

\end{tikzpicture}
    \caption{Starting from the left, the neighbor graph for the dendrite of Figure \ref{fig:attractorvsJulia}, the Sierpi\'{n}ski triangle, and the Hexagasket. The position of the vertices reflects the position of the $h(A_\lambda)$ relative to $A_\lambda$.}
    \label{fig:finitelyramified}
\end{figure}
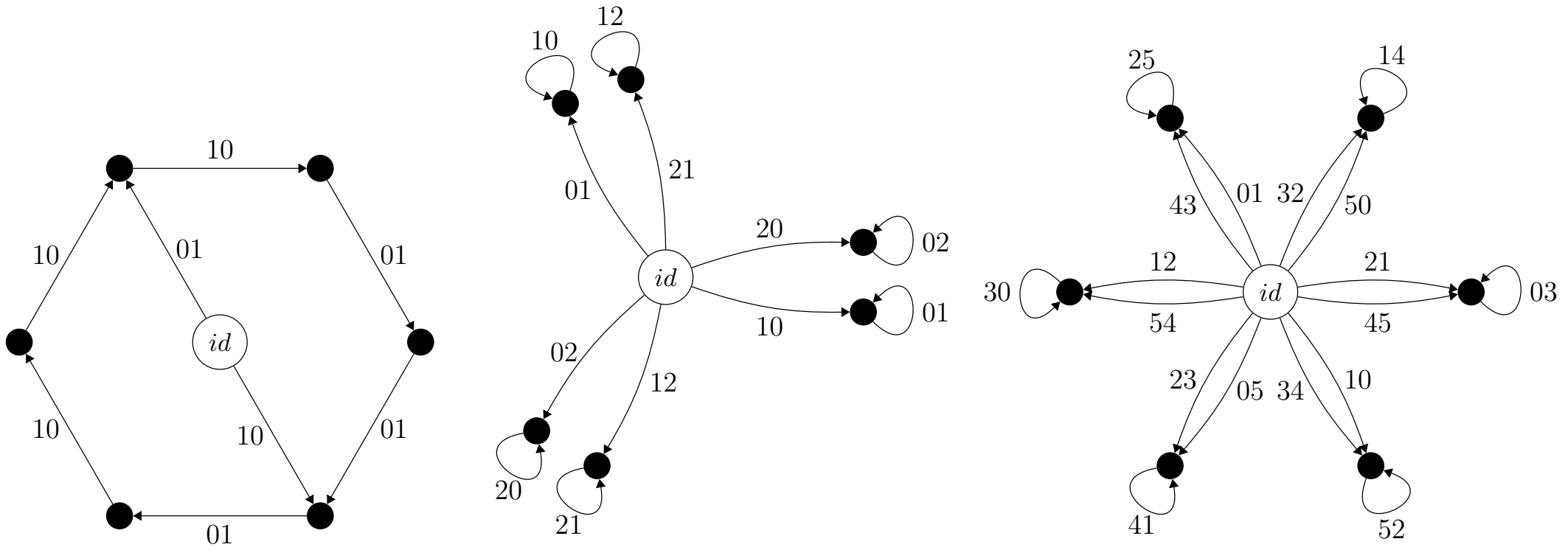

For example, in the neighbor graph of the Sierpi\'{n}ski triangle shown in Figure \ref{fig:finitelyramified}, the infinite path $01,10,10,\ldots$ corresponds to the point $p=1+\xi_3$ in the attractor since $\lambda=0.5$:\[\pi_\lambda(\mathbf{0}\overline{\mathbf{1}})=1+\frac{\lambda \xi_3}{1-\lambda}=\xi_3+\frac{\lambda}{1-\lambda}=\pi_\lambda(\mathbf{1}\overline{\mathbf{0}})\] and since the two addresses differ at the first symbol, $p$ lies in the intersection of $A_0$ and $A_1$. Another infinite path $10,01,01,\ldots$ through a different vertex in the graph gives the same equation $\pi_\lambda(\mathbf{1}\overline{\mathbf{0}})=\pi_\lambda(\mathbf{0}\overline{\mathbf{1}})$. Since these two paths are the only ones involving the symbols $0$ and $1$, $p$ is the only element of the intersection $A_0\cap A_1$. 

In the Hexagasket (also shown in Figure \ref{fig:finitelyramified}), there are pairs of infinite paths that go through the same vertex, but they represent different points in the attractor. For example, the infinite paths $21,03,03,\cdots$ and $45,03,03,\cdots$ correspond to points in $A_1\cap A_2$ and $A_4\cap A_5$, respectively.

The number of vertices in the neighbor graph indicates a different property of the attractor. Except the root, each vertex represents a \emph{boundary set} of $A_\lambda$, and the first labels $s_1s_2\cdots$ of the infinite path starting in $h$ give the addresses of the points of this boundary sets. Consider a the edge labeled $s_1t_1$ from the root to the vertex $h$, then \[D=A_{s_1}\cap A_{t_1}=f_{s_1}(A_\lambda)\cap f_{t_1}(A_\lambda)\neq\emptyset \iff f_{s_1}^{-1}(D)=A_\lambda\cap h(A_\lambda)\neq\emptyset.\]
\begin{thm}[\cite{BandtMesing2009}]
    \label{thm:dynamicalboundary}
    For a finite type fractal with open set condition, the addresses of the boundary set corresponding to $h$ form the regular language $L_h$ defined by the first labels on $NG_\lambda\setminus\{id\}$ with intitial vertex $h$.
\end{thm}
In our examples, there are two distinct neighbors $f_{2}^{-1}f_{0}(A_\lambda)$ and $f_{1}^{-1}f_{0}(A_\lambda)$ of the Sierpi\'{n}ski triangle $\Phi_3(A_\lambda)$ meeting at the boundary set $\{\pi_\lambda(\overline{\mathbf{0}})\}$; while the Hexagasket $\Phi_6(A_\lambda)$ has exactly one neighbor $f_{2}^{-1}f_{1}(A_\lambda)=f_{4}^{-1}f_{5}(A_\lambda)$ at the boundary set $\{\pi_\lambda(\overline{\mathbf{0}})\}$.

The union of all the boundary sets was called \emph{dynamical boundary} by Mor\'{a}n \cite{Moran1999} and he proved that if $\Phi_N(A_\lambda)$ has non-empty interior and satisfies the open set condition, then its dynamical boundary is its topological boundary. Sirvent \cite{Sirvent2014} in the study of some particular type of IFS, which include $\{f_j\circ s_{\sigma(j)}\}_{j=0}^1$, gave a criterion for certain points to be in the dynamical boundary. In Section \ref{sect:IIFSAndConjugacy} we reinterpret the dynamical boundary as the post-critical set of the induced expanding map $q_\lambda$ on $A_\lambda$.

Bandt and Mesing distinguish three classes of finite type fractals: finitely
ramified, infinitely ramified and overlapping fractals depending on how the vertices are connected. We are interested only in the first one, but we need some more terminology: a \emph{terminal cycle} in $NG_\lambda$ is a directed cycle without any diagonal edge or double edge or any edge leading out of the cycle. A vertex $h$ of $NG_\lambda$ is called \emph{terminal} if only terminal cycles can be reached from it. A finite type fractal is called \emph{finitely ramified} if it satisfies $\#(A_j\cap A_k)<\infty$ for every $j,k\in\Sigma_N$.
\begin{thm}[\cite{BandtMesing2009}]
\label{thm:finitelyramified}
    A finite type fractal is finitely ramified if and only if all cycles in the neighbor graph are terminal cycles.
\end{thm}
The Sierpi\'{n}ski triangle and the Hexagasket are examples of finitely ramified fractals. In Figure \ref{fig:finitelyramified} there is also the neighbor graph associated to the attractor $\Phi_2(A_\lambda)$ when $\lambda\approx 0.366+0.520i$, which is also finitely ramified.

\section{Complex Dyamics}
\label{sect:ComplexDynamics}
The main object of study in complex dynamics is the behavior of a point $z$ under iteration of an holomorphic map $F$. The {\em orbit} of $z \in \mathbb{C}$ is the set of iterates $\{z,\,F(z),\,F^2(z),\,\ldots\}$. Some orbits may repeat in a loop, some may wander around the plane, some others may diverge to infinity. Define the {\em filled Julia set} as
\[K(F) := \{ z\in\mathbb{C}~\mid~\text{the orbit of } z \text{ under } F \text{ is bounded}\,\}\]
and the \emph{Julia set} $J(F)$ is simply the boundary of $K(F)$. As for the IFS setting, the connectedness of $J(F)$ is a key property to be studied. For this reason we define the connectedness locus of the function \[\mathcal{M}(F)=\{z\in\mathbb{C}~\mid~J(F)\text{ is connected}\}\] commonly referred to as the Mandelbrot set of $F$. Note that in the future we may drop the dependence of $F$ in the notation of these sets if it is clear from context.

\subsection{The Mandelbrot set and Misiurewicz parameters}
\label{subsect:mandelbrotsetsforQuadraticmaps}
When the holomorphic function $F$ is the quadratic map $Q_c(z)=z^2+c$ with $c\in\mathbb{C}$, there is a famous dichotomy which motivated our analogy with the IFS setting when $N=2$:
\begin{thm}[\cite{Orsay}]
For $c\in\mathbb{C}$, either 
\begin{enumerate}
    \item $0 \in K(Q_c)$, and $K(Q_c)$ is connected or
    \item $0 \not\in K(Q_c)$, and $K(Q_c)$ is a Cantor set.
\end{enumerate}
\end{thm}
This results in the following characterization of the Mandelbrot set
\[\mathcal{M}(Q_c) := \{ z\in\mathbb{C}~\mid~ 0 \in K(Q_c) \},\]
which also makes possible a computable criterion to create pictures of it.
We already mentioned some properties of $\mathcal{M}(Q_c)$:
\begin{itemize}
    \item[-] $\mathcal{M}(Q_c)$ is connected. The crowning result of the celebrated Orsay notes of Douady and Hubbard \cite{Orsay}.
    \item[-] MLC Conjecture. A most sought open problem that states that the Mandelbrot set is locally connected. 
    \item[-] Both $\mathcal{M}(Q_c)$ and $K(Q_c)$ are full in the sense that their complements are connected.
    \item[-] A parameter $c \in \mathcal{M}(Q_c)$ is {\em Misiurewicz} if the critical orbit terminates in a periodic cycle that does not contain $0$ itself; i.e., is strictly {\em pre-periodic}. The corresponding Julia set is locally connected but contains no simple closed curves; that is $J(Q_c)$ is a dendrite.
    \item[-] \cite{Lei1990} If $c \in \mathcal{M}(Q_c)$ is Misiurewicz, then $\mathcal{M}(Q_c)$ ``looks locally like'' the filled Julia set $K(Q_c)$: there are two sequences $\{\varphi_n\},\,\{\psi_n\}$ of homotheties centered at $c$ so that $\overline{\mathbb{D}}_c \cap \varphi_n(\mathcal{M})$ and $\overline{\mathbb{D}}_c\cap \psi_n(K)$ are asymptotically close in the Hausdorff metric, where $\mathbb{D}_c$ is the unit disk centered at $c$.
  \end{itemize}
Misiurewicz parameters are dense in $\partial\mathcal{M}(Q_c)$ \cite{McMullen1988,Lei1990}, and for such $c$:
\begin{itemize}
    \item[-] every cycle of $Q_c$ is repelling \cite{Orsay};
    \item[-] $J(Q_c)$ is connected and locally connected \cite{Orsay};
    \item[-] $J(Q_c)$ has no parabolic cycles \cite{Orsay};
    \item[-] $Q_c$ is semi-hyperbolic \cite{CarlesonJonesYoccoz1994}.
\end{itemize}


\subsection{External angles, kneading sequences, and core entropy}
\label{subsect:raysKSentropy}
Misiurewicz parameters admit also a description via external angle. According to Duoady and Hubbard \cite{Orsay}, there is a unique analytic isomorphism $\widehat\phi$ mapping the exterior of the closed unit disk to the exterior of the Mandelbrot set. The \emph{external ray} $R_{\mathcal{M}}(\theta_0)$ with angle $\theta_0$ is the $\widehat\phi$-image of $\theta=\theta_0$. If $\theta\in\mathbb{R}/2\pi\mathbb{Z}$ is rational, the external ray $R_{\mathcal{M}}(\theta)$ lands at a parameter $c(\theta)$. That is $\lim_{r\to 1^+}\widehat\phi(r\exp(2\pi i \theta))$ exists and it is a unique point on $\partial\mathcal{M}$. The parameter $c(\theta)$ is Misiurewicz precisely when $\theta$ is strictly preperiodic under the doubling map $\theta\mapsto2\theta\bmod1$. It follows that, Misiurewicz parameters are exactly those landed on by angles
\[\theta=\frac{t}{2^k(2^n-1)},\qquad n,k\geq1\text{ and }t\text{ odd}.\]

There exist external rays in the dynamical plane as well, see \cite{Milnor2006}. In fact, when $J_c:=J(Q_c)$ is locally connected every external ray lands at a point in $J_c$. Suppose that $c\in J_c$ is the landing point of an external ray with angle $\theta$. Then the rays with angle $\theta/2$ and $(\theta+1)/2$ land at $0$ and partition $J_c$ into three pieces $\{\{0\},J^{\mathbf{1}},J^{\mathbf{0}}\}$ such that $c\in J^{\mathbf{1}}$. The itinerary of the critical value with respect to this partition is the \emph{kneading sequence} of $c$. Such sequence coincides with the itinerary of $\theta$ under the angle doubling map on $\mathbb{R}/2\pi\mathbb{Z}$ with respect to the partition $\{\mathbf{1}:=(\theta/2,(\theta+1)/2),\mathbf{0}:=((\theta+1)/2,\theta/2)\}$.
This combinatorial invariant is a key ingredient in the determining the topological conjugacy between $(A_\lambda,q_\lambda)$ and $(J_c,Q_c)$. 

Indeed, when $\lambda\in\mathcal{T}_2$ the attractor $A_\lambda$ can be partitioned in three pieces $\{\{0\}, A_1\setminus\{0\}, A_0\setminus\{0\}\}$. The map $q_\lambda:A_\lambda\to A_\lambda$ is always defined (see Section \ref{sect:IIFSAndConjugacy}) for such parameters, therefore we can define the kneading sequence for $q_\lambda$ as the itinerary of $q_\lambda(0)$ with respect to this partition. If such sequence equals the kneading sequence for the angle doubling map for some angle $\theta$, then it is called \emph{admissible}.
\begin{thm}[\cite{ErogluRohdeSolomyak2010}] 
Suppose $\lambda\in\mathcal{T}_2$ is such that the following hold.
\begin{itemize}
    \item[-] The kneading sequence for $q_\lambda$ is admissible and corresponds to an angle $\theta$.
    \item[-] There is a quadratic Julia set $J_c$ which is locally connected with external angle $\theta$. 
\end{itemize}
Then $(A_\lambda,q_\lambda)$ and $(J_c,Q_c)$ are conjugate.
\end{thm}
This results hints to the following open question:
\begin{quest}
    What is the set of parameters $c$ in the quadratic Mandelbrot set for which there exists a corresponding $\lambda\in\mathcal{T}_2$?
\end{quest}
The point $c=i$ is an example of a Misiurewicz parameter which does not correspond to a $\lambda\in\mathcal{T}_2$. The orbit of $0$ lands on the $2$-cycle $\{i-1,-i\}$ after one iteration. The kneading sequence is then $\mathbf{1}\overline{\mathbf{10}}$ corresponding to the angle $1/6$. If there were $\lambda\in\mathcal{T}_2$ with the same kneading sequence, according to \cite{ErogluRohdeSolomyak2010}, $\lambda$ would equal $(1+\sqrt{7}i)/4$. However, this parameter is not in $\mathcal{T}_2$ since its attractor has non-empty interior and it is known to tile the plane, that is $A_\lambda$ is not a dendrite.  Other examples of Misiurewicz parameter which have no corresponding $\lambda\in \mathcal{T}_2$ are those associated to external angles $1/8$, $1/10$, and $3/10$. In fact, repeating the argument above we obtain parameters $\lambda$ with modulus greater than or equal to $1/\sqrt{2}$ and hence attractors with substantial overlap.

For parameters $c$ whose orbit of $0$ is finite, there is a finite and forward-invariant topological tree embedded in $K(Q_c)$ called the \emph{Hubbard tree} $H_c$ \cite{Orsay}. The action of $Q_c$ on $H_c$ determines the topological conjugacy class of $Q_c$. In reality Hubbard trees can be defined more abstractly and not all Hubbard trees are realized by quadratic polynomials. Thurston \cite{Thurston2014} defined the \emph{core entropy}, $h(c)$ as the the topological entropy of $Q_{c}$ restricted to its Hubbard tree $H_{c}$. He conjectured that $h(c):\mathcal{M}\to[0,\log2]$ is a continuous function of $c$. Tiozzo \cite{Tiozzo2016} proved continuity with respect to the external angle $\theta$ and Dudko and Schleicher \cite{DudkoSchleicher2020} settled it. In fact, they also proved
\begin{thm}[\cite{DudkoSchleicher2020}]
The core entropy function $h:\mathbb{R}/2\pi\mathbb{Z}\to[0,\log2]$ has the following properties.
\begin{itemize}
    \item[-] Every dyadic angle is an isolated local maximum of the entropy function.
    \item[-] Conversely, every local maximum of $h$ is dyadic.
    \item[-] Within every wake, the entropy function has a unique global maximum, and it occurs at the unique dyadic angle of lowest denominator in the wake.
\end{itemize}    
\end{thm}
A \emph{wake} is a closed interval in $\mathbb{R}/2\pi\mathbb{Z}$ bounded by two rational angles so that the corresponding rays land together.

Remarkably, the system $(A_\lambda,q_\lambda)$ for $\lambda\in\mathcal{T}_2$ satisfy all the axioms of an abstract Julia set of Bruin and Schleicher\cite{BruinSchleicher2008}. Moreover, the connected hull of the orbit of $0$ under $q_\lambda$ is exactly the Hubbard tree. The work in \cite{ErogluRohdeSolomyak2010} also show that the dynamics of $q_\lambda$ cyclically permutes the `arms' of the Hubbard tree at periodic branch points, a necessary and sufficient condition on Hubbard trees to be realizable by a quadratic polynomial.

More independent experimental data reveals that the Misiurewicz parameters $c$ for which a corresponding $\lambda\in\mathcal{T}_2$  seems to exist have relatively high core entropies within their limb. For an angle $\theta$ in the wake of a $p/q$-limb, let $\theta^*$ be the dyadic angle of lowest denominator in that wake. In the $1/3$-limb, bounded by the angles $1/7$ and $2/7$, we have $\theta^*=1/4$ and all the angles in our sample for which a correspondence was found satisfy $h(\theta)\geq 0.8 h(1/4)$. In the $1/2$-limb, bounded by the angles $1/3$ and $2/3$, $\theta^*=1/2$ and they satisfy $h(\theta)\geq 0.68 h(1/2)$. Unfortunately, these bounds are far from sufficient since they are also satisfied by 82\% and 45\%, respectively, of the angles for which no correspondence was found. 
\begin{figure}[!htb]
    \centering
    \includegraphics[scale=0.65]{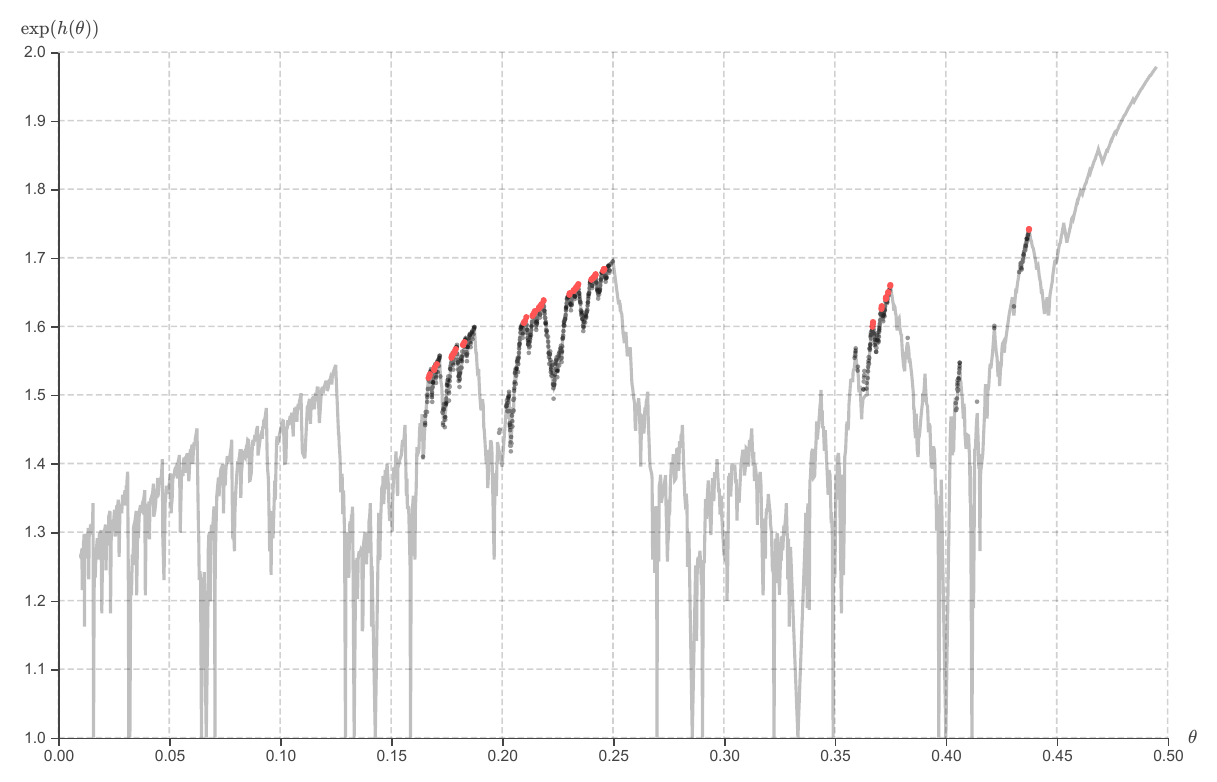}
    \caption{In gray the exponential of the core entropy $h(\theta)$  for $0<\theta<1/2$ as studied by Tiozzo. Marked in black is $\exp(h(\theta))$ for a large class of $\theta$ associated to Misiurewicz parameters. The highlighted one correspond to those parameters for which there seems to be a corresponding $\lambda\in\mathcal{T}_2$.}
    \label{fig:entropy}
\end{figure}

The plot in Figure \ref{fig:entropy} was obtained in the following way: we considered the kneading sequences for Misiurewicz parameters $c(\theta)$ where $\theta=t2^{-k}(2^n-1)^{-1}$ with odd $t$ and specified $k\geq1$ and $n\geq0$, recovered the possible corresponding $\lambda\in\mathbb{B}_2$ (in the first quadrant) using the method described in \cite{ErogluRohdeSolomyak2010} and marked only those with modulus strictly less than $1/\sqrt{2}$. For each candidate we then computed the neighbor graph (Section \ref{subsect:neighborgraphs}) until an iteration produced no new vertices nor edges or it reached a maximum depth. Those that reached maximum depth were discarded. We labeled ``working'' only those $\lambda$ for which there is no path in the resulting $NG_\lambda$ leading back to the identity (Theorem \ref{thm:OSC}) and no edge with label $00$ or $11$ survives the pruning (Lemma \ref{lem:no0fordendrite}).

One final comment about the experimental data: in accordance with a remark in \cite[Section~5.4]{ErogluRohdeSolomyak2010}, the kneading sequence $\nu=v_1v_2\cdots$ of every ``working'' angle we have found has no two consecutive $\mathbf{1}$s after $v_2$.
\begin{conj}
    Suppose the two dynamical systems $(A_\lambda,q_\lambda)$ and $(J_c,Q_c)$ are (quasisymmetrically) conjugate, with $\operatorname{Re}\lambda>0$ and $\operatorname{Im}\lambda\geq0$. Then the kneading sequence of $c$ will start with $\mathbf{11}$ and $\mathbf{11}$ will never occur again.
\end{conj}


\subsection{The parameter plane for McMullen maps}
\label{subsect:mandelbrotsetsforMcMullenmaps}
Fix $n\geq2$, $m\geq1$ and let $N=n+m$. The rational map $R_{n,m,c}(z)=z^n+cz^{-m}$ for $c\in\mathbb{C}$ has a superattracting fixed point at $\infty$, a pole of order $m$ at $0$, and $N$ critical points at \[p_j=\sqrt[N]{\lvert c\rvert\frac{m}{n}}\exp((\arg(c)+2\pi ij))/N),~0\leq j\leq N-1.\] Observe that if we conjugate the McMullen maps with $z\mapsto c^{1/N}z$, then for the resulting family $F_{n,m,\mu}(z):=\mu R_{n,m,1}(z)$ with $\mu=c^{(n-1)/N}$, the critical points will lie on the circle of radius $\lvert z\rvert=(m/n)^{1/N}$ and are independent from $c$. Due to the symmetries below, Steinmetz \cite{Steinmetz2006} prefers to study $R_{n,m,c^N}$ and $R_{n,m,c^{N(n-1)}}$ rather than $R_{n,m,c}$, but we will follow Devaney's lead.

There are three basic symmetries to be aware:
\begin{itemize}
    \item[\emph{Symmetry 1.}] The map $R_{n,m,c}$ is conjugate to $R_{n,m,\overline{c}}$ via $z\mapsto\overline{z}$.
    \item[\emph{Symmetry 2.}] If $\xi_N$ is a $N$th root of unity, then $R_{n,m,c}(\xi_Nz)=\xi_N^nR_{n,m,c}(z)$.
    \item[\emph{Symmetry 3.}] Let $\eta=\xi_N^{1/(n-1)}$ and $\nu=\eta^{N}$. Then for all $k\geq1$\[R_{n,m,\nu c}^k(\eta z)=\eta^{n^k}R_{n,m,c}^k(z).\]
\end{itemize}
The first symmetry implies that the parameter plane is symmetric with respect to complex conjugation and the Julia sets $J(R_{n,m,c})$ and $J(R_{n,m,\overline{c}})$ are homeomorphic. The second one implies that the Julia set $J(R_{n,m,c})$ is symmetric under rotation by $\xi_N$ and consequently the critical orbits behave symmetrically with respect to this rotation. 

The last symmetry allows to determine the orbit diagram of $R_{n,m,\nu c}$ from the one of $R_{n,m,c}$. In particular, if $p$ is a critical point for $R_{n,m,c}$ then $\eta p$ is one of $R_{n,m,\nu c}$. Therefore, the orbits of the critical points of $R_{n,m,c}$ and $R_{n,m,\nu c}$ behave symmetrically with respect to rotation by some power of $\eta$. \emph{Symmetry 3} also implies that the full basin of $\infty$ for $R_{n,m,c}$ is homeomorphic to that of $R_{n,m,\nu c}$ under rotation by $\eta$. Since this basin is completely invariant, the Julia sets of one map is the rotation of the other under $z\mapsto\eta z$. 

Unlike the quadratic family $Q_c$, there are three ways the orbit of $0$ can tend to $\infty$ which generate different topologies for $J_c:=J(R_{n,m,c})$. Since the point at $\infty$ is superattracting fixed point, it is surrounded by an immediate basin of attraction, which is denoted by $B_c$. Now the full basin of attraction may consist of infinitely many disjoint preimages of $B_c$. In particular, the component $T_c$ of the full basin containing $0$ may be disjoint from $B_c$. If $T_c\cap B_c=\emptyset$, then $R_{n,m,c}$ is $n$-to-$1$ on $B_c$ and $m$-to-$1$ on $T_c$ and thus there are exactly two preimages of $B$, namely $B_c$ and $T_c$. It follows that $T_c$ must have disjoint preimages under $R_{n,m,c}^k$ for $k\geq1$, so the full basin of $\infty$ must have infinitely many components.

\begin{thm}[The escape trichotomy \cite{DevaneyLookUminsky2005}]
    Suppose the orbits of the critical points tend to $\infty$. 
    \begin{enumerate}
        \item if $v_0$ lies in $B_c$, then $J_c$ is a Cantor set;
        \item if $v_0$ lies in $T_c$, then $J_c$ is a Cantor set of concentric simple closed curves, each one of which surrounds the origin;
        \item in all other cases, $J_c$ is a connected set, and if $R_{n,m,c}^k(v_0)\in T_c$ where $k\geq1$, then $J_c$ is a Sierpi\'{n}ski curve.
    \end{enumerate}
\end{thm}
Our interest however is for when the orbit of $p_0$ remains bounded. In this case, the Julia set is connected \cite{DevaneyRussell2013}. The set of parameters for which the critical orbits either do not escape or else enter the trap door $T_c$ at iteration $2$ or higher is the connectedness locus $\mathcal{M}(R_{n,m,c})$. A complete topological description of Julia sets of a general McMullen family $R_{n,n,c}(z)+a$ for $a\in\mathbb{C}$ was given in \cite{XiaoQiuYin2014}. 

We further restrict our attention to those values of $c$ for which the following conditions holds:
\begin{itemize}
    \item[-] each critical point lies on the boundary of $B_c$
    \item[-] each of the critical points are strictly preperiodic.
\end{itemize}
These parameters are called Misiurewicz. In general, Misiurewicz parameters form a dense subset in the unstable locus of families of rational-like maps \cite{Tan1998}. Devaney, Rocha, and Siegmund \cite{DevaneyRochaSiegmund2007} showed that the set of Misiurewicz parameters is a dense subset of the boundary of the connectedness locus $\mathcal{M}(R_{n,m,c})$. Moreover, the Julia set they proved that for Misiurewicz parameters the Julia set is a generalized Sierpi\'{n}ski gasket.

We use the definition of generalized Sierpi\'{n}ski gasket given in \cite{HeXiaoYang2025} since it sharpens and unifies the original Devaney--Rocha--Siegmund's description \cite{DevaneyRochaSiegmund2007} with the $N=1$ case of Luo and Ntalampekos \cite{LuoNtalampekos2026}. 
\begin{defn}
    \label{def:generalizedcarpet}
    A locally connected continuum $K\subset\mathbb{C}$ with empty interior is a \emph{generalized Sierpi\'{n}ski gasket}, if 
    \begin{enumerate}
        \item each complementary component of $K$ is a Jordan domain;
        \item there exists $N\geq1$ such that the boundaries of any two complementary components of $K$ share at most $N$ points;
        \item no three complementary components have a common boundary point;
        \item the contact graph of $K$, obtained by assigning a vertex to each complementary component and an edge if two touch, is connected.
    \end{enumerate}
\end{defn}
We remark that this is a purely topological definition and not based on dynamics. It says nothing about which branched self-cover structure sits on top of it. For example, two Misiurewicz parameters can produce Julia sets that are homeomorphic as \emph{sets} while carrying genuinely non-conjugate dynamics.

\subsubsection{Parameter rays} As in the quadratic case, the boundary of the escape locus $H_0\subset\mathbb{C}\cup\{\infty\}$, i.e. the complement of $\mathcal{M}$, admits a theory of parameter rays. For the symmetric subfamily $R_{n,n,c}$, Qiu, Roesch, Wang, and Yin \cite{QiuRoeschWangYin2015} construct a canonical parametrization $\nu:\mathbb{R}/2\pi\mathbb{Z}\to\partial H_0$, with $\nu(\theta)$ defined as the landing point of the parameter ray $\mathcal{R}_\mathcal{M}(\theta)$, and prove that $\nu(\theta)\in\partial H_0$ has $\partial B$, the boundary of the immediate basin of $\infty$, containing either the critical set or a parabolic cycle, with the latter occurring exactly when $\theta$ is periodic under $\tau(\theta)=n\theta\bmod1$.
 
In \cite{DevaneyRochaSiegmund2007} is established that $J(R_{n,m,c})$ is a generalized Sierpi\'{n}ski gasket only for the critically finite, Misiurewicz parameters. The authors of \cite{QiuRoeschWangYin2015} had shown $c\in\partial H_0$ if and only if $\partial B_c$ contains either the critical set or a parabolic cycle. He, Xiao, and Yang \cite{HeXiaoYang2025} prove that the gasket conclusion, in fact, holds for every $c\in\partial H_0$, without any finiteness assumption on the critical orbit:
 
\begin{thm}[\cite{HeXiaoYang2025}, Theorem 1.1]
\label{thm:hxy}
If $c\in\partial H_0$, then the Julia set of $R_{n,n,c}$ with $n\geq3$ is a generalized Sierpi\'{n}ski gasket.
\end{thm}
 
In particular this covers parameters with a parabolic cycle on $\partial B_c$ (the cusps of $\nu$), as well as critically infinite parameters such as those with a Siegel disk or Cremer point, going well beyond the strictly preperiodic Misiurewicz parameters that are our main concern in this survey. The authors explicitly state that they expect Theorem \ref{thm:hxy} to extend to the asymmetric family $R_{n,m,c}$, $n\neq m$, but note that their proof relies on the machinery of \cite{QiuWangYin2012,QiuRoeschWangYin2015}, which is developed only for the symmetric case $n=m$.

The critical points of $R_{n,m,c}$ solve $z^{n+m}=mc/n$, giving $n+m$ simple critical points regardless of whether $n=m$. The local degrees at two distinguished points is what differs the two cases: $\infty$ has local degree $n$, while $0$ is a pole of order $m$. When $n=m$ these degrees coincide, which is what allows an involution intertwining the basin of $\infty$, $B_c$ and the trap door, $T_c$; when $n\neq m$ such involution can no longer exist. On the other hand, the escape locus $H_0$ is governed entirely by the local degree at $\infty$: the B\"{o}ttcher coordinate near $\infty$ satisfies $\widehat\phi(R_{n,m,c}(z))=\widehat\phi(z)^n$ independent of the pole order $m$. This suggests that the cusp and Misiurewicz dichotomy of \cite{QiuRoeschWangYin2015}, and consequently the gasket conclusion of Theorem \ref{thm:hxy}, should extend to the asymmetric family, but for now this remains an open problem.


\subsection{The model of a generalized gasket}
\label{subsect:gasketmodel}
For Misiurewicz parameters producing a generalized Sierpi\'{n}ski gasket, Devaney, Moreno Rocha, and Siegmund \cite{DevaneyRochaSiegmund2007} construct an explicit combinatorial model of which we recall the construction.
 
Recall that $B_c$ indicates the immediate basin of $\infty$ and $T_c$ the trap door, the component of $R_{n,m,c}^{-1}(B_c)$ containing $0$; write $\beta=\partial B_c$ and $\tau=\partial T_c$. When $c$ is Misiurewicz, both $\beta$ and $\tau$ are simple closed curves, and $\beta\cap\tau$ consists precisely of the $n+m$ critical points, called the corners of the trap door \cite{DevaneyRochaSiegmund2007}. Together with the strict preperiodicity of the critical points, this is how \cite{DevaneyRochaSiegmund2007} define 
Misiurewicz parameters.
 
Writing $\tau^k$ for the union of all components of $R_{n,m,c}^{-k}(\tau)$, each component of $\tau^k$ is again a simple closed curve with $n+m$ corners, and the nested sets
\[K_0 = \widehat{\mathbb{C}}\setminus B_c,\qquad K_{j+1} = K_j \setminus R_{n,m,c}^{-j}(T_c)\]
satisfy $J_c=\bigcap_{j\geq0} K_j$, exhibiting $J_c$ directly as a generalized Sierpi\'{n}ski gasket. The precise combinatorial rule governing how these polygons nest, and the object that will matter most for our purposes, is given in the following lemma.
\begin{lem}[\cite{DevaneyRochaSiegmund2007}, Lemma 5.1]
\label{lem:mcorners}
Let $C$ be a component of $\tau^k$, $k\geq1$. Then exactly $m$ of its $n+m$ corners lie on an edge of a single component of $\tau^{k-1}$.
\end{lem}
 
Applying Lemma \ref{lem:mcorners} produces an explicit geometric model. Starting from the unit circle with the map $z\mapsto z^{n}$, to which $R_{n,m,c}$ is conjugate near $\beta$, one successively ``straightens'' $\tau,\tau^1,\tau^2,\dots$ into nested homeomorphic images built purely from the itinerary of the critical orbit, producing a homeomorphism between $J_c$ and this model. Devaney, Moreno Rocha, and Siegmund further show that two Misiurewicz Julia sets in the same family $R_{n,m,c}$ are homeomorphic if and only if their models coincide up to complex conjugation.

\section{Invertible IFS and dynamical conjugacy}
\label{sect:IIFSAndConjugacy}
Recall that each attractor $A$ of an IFS with $N$-maps $\{f_j\}_{j=0}^{N-1}$ has associated a symbol space $\Sigma_N=\{0,1,\cdots,N-1\}$ and a continuous projection map $\pi:\Sigma_N^\infty\to A$.

\begin{defn}[\cite{Kameyama1993}]
A self similar set $A$, attractor of an IFS $\{f_j\}_{j=0}^{N-1}$, is called \emph{post-critically finite} (p.c.f. for short) if 
\begin{enumerate}
    \item All $f_j$ are injective.
    \item For $i\neq j$, $\#(A_i\cap A_j)\leq1$.
    \item For points $c\in\mathcal{O}=\cup_{i\neq j}(A_i\cap A_j)$ the addresses in $\pi^{-1}(c)$ are preperiodic by the right shift map on the symbol space $\Sigma_N$.
\end{enumerate}
\end{defn}
We remark that the finitely ramified fractal defined in Section \ref{subsect:neighborgraphs} are p.c.f. fractals if all cycles in the neighbor graph are terminal cycles.

\begin{defn}[\cite{Kameyama1993}]
The attractor $A$ of an IFS $\{f_j\}_{j=0}^{N-1}$, is called \emph{invertible post-critically finite} (abbreviated as i.p.c.f.) if the following conditions are satisfied.
\begin{enumerate}
    \item All $f_j$ are injective
    \item For $i\neq j$, $\#(A_i\cap A_j)<\infty$.
    \item There exists a map $q$ on $A$ to itself such that $q\circ\pi=\pi\circ\tau$, where $\tau$ is the right-shift map on $\Sigma_N$.
    \item All points in $\mathcal{O}=\bigcup_{i\neq j}(A_i\cap A_j)$ are preperiodic by $q$.
\end{enumerate}
\end{defn}
Note that the third condition implies that $q=f_j^{-1}$ on each $A_j$. The set $\mathcal{O}$ is the critical set of $q$ in the sense that it fails to be a local homeomorphism precisely there. By a result of Bandt and Rao \cite{BandtRao2007}, if $A_i\cap A_j$ is a finite set for $i\neq j$, then the open set condition holds which, in turn, implies that the points in $\mathcal{O}$ do not have a recurrent address. The IFS being \emph{critically non-recurrent}, i.e. that for every $p\in\mathcal{O}$ the post-critical set $\{q^k(p)\}_{k\ge1}$ does not contain $p$, is weaker than the last condition in the definition of i.p.c.f.. In symbolic dynamics terms the attractor $A$ is called a \emph{$N$-factor}, $\pi_\lambda$ a \emph{factor map}, and $A=\Sigma_N^\infty/\sim$ is $\tau$-invariant since $\mathbf{s}\sim\mathbf{t}\iff\tau(\mathbf{s})\sim\tau(\mathbf{t})$, see \cite{BandtKeller1991}.
\begin{thm}[\cite{Kameyama1993}]
\label{thm:kameyamaPolynomial}
    Let $A$ be the attractor of an IFS of similarity maps and suppose it is simply connected and i.p.c.f. with map $q$. Then there exists a polynomial map $g$ with Julia set $J(g)$ and a homeomorphism $h:J(g)\to A$ conjugating $g$ and $q^n$ for some integer $n$.
\end{thm}
In \cite{Kameyama2000}, Kameyama extended the result to attractors that are not simply connected:
\begin{thm}[\cite{Kameyama2000}]
\label{thm:kameyamaRational}
    Let $A$ be the attractor of an IFS of similarity maps and suppose it is i.p.c.f. with map $q$. Then there exists a postcritically finite rational map $g$ with Julia set $J(g)$ and a homeomorphism $h:J(g)\to A$ conjugating $g$ and $q$. Moreover, the homeomorphism $h$ can be extended to the whole Riemann sphere. 
\end{thm}

Ero\u{g}lu, Rohde, and Solomyak strengthen Kameyama's result in particular cases: the homeomorphism is quasisymmetric if the attractor $A$ is of bounded turning and the rational map $g$ is \emph{semi-hyperbolic}, i.e. that it does not have parabolic periodic points and that the critical points are not recurrent. 

The map $q$ is only defined a priori on the i.p.c.f. attractor $A$. Kameyama's proof of Theorem \ref{thm:kameyamaPolynomial} passes through Thurston's topological characterization of p.c.f. rational maps, proven in \cite{DouadyHubbard1993}. A postcritically finite branched covering of the sphere with hyperbolic orbifold is Thurston equivalent to a rational map if and only if it admits no Thurston obstruction, and in that case the rational map is unique up to a M\"{o}bius conjugation. Kameyama's contribution is to use the combinatorial structure defined by the identifications in $\Sigma_N^\infty$ to extend $q$ to the whole plane. By Theorem \ref{thm:equivalentaddresses}, for attractors of finite type the combinatorial structure is encoded in the neighbor graph.

The upgrade by Ero\u{g}lu--Rohde--Solomyak to quasisymmetry of the conjugacy follows the strategy of McMullen and Sullivan \cite{McMullenSullivan1998} for hyperbolic maps, adapted to a non-hyperbolic setting by Przytycki and Rohde \cite{PrzytyckiRohde1999}. The proof shows that the conjugacy maps balls to uniformly roundish sets\footnote{A set $S\subset X$ is $C$-roundish at $x$ if $\operatorname{diam}(S)\leq C\inf\{\lvert x-y\rvert~\mid~y\in X\setminus S\}$.}. In particular, a small ball is squeezed between two pullback components whose size is controlled on the Julia side by semi-hyperbolicity of $g$ \cite{CarlesonJonesYoccoz1994} and on the IFS side by the fact that $q$ is piecewise a similarity with non-recurrent overlap set. Quasisymmetric maps preserve bounded turning, and since $J(g)$ has it by \cite{CarlesonJonesYoccoz1994}, $A$ must also have it for any quasisymmetric conjugacy to exist. In fact, the bounded turning property of $A$ controls the pullback components on the IFS side. In the quadratic case, the conjugacy extends to a quasiconformal homeomorphism of the sphere \cite[Theorem 1.1]{ErogluRohdeSolomyak2010}, because the complement of the attractor is a John domain \cite{NakkiVaisala1991}.

Ha\"{i}ssinsky and Pilgrim \cite{HaissinskiPilgrim2012} obtain the quasisymmetry statement from the principle of the conformal elevator. They show that the invertible IFS $\Phi_{2,\sigma}$, with the Euclidean metric, are metrically coarse expanding conformal (CXC) dynamical systems, and so are semi-hyperbolic rational maps on their Julia sets.

An alternative route to the same kind of conjugacy statement obtained via Ero\u{g}lu--Rohde--Solomyak and Ha\"{i}ssinsky--Pilgrim, can be found in the work of Belk and Forrest. In \cite{BelkForrest2025}, the authors develop an independent theory of quasisymmetries for finitely ramified fractals, a notion introduced by Teplyaev \cite{Teplyaev2008} as a common generalization of Kigami's postcritically finite self-similar sets \cite{Kigami2001}. Their approach constructs, on such fractals, a distinguished quasisymmetric equivalence class of undistorted metrics, and shows that piecewise defined homeomorphisms respecting the cell structure are automatically quasisymmetries.

\subsubsection{Known examples}
\label{subsubsect:knownexamples}
In what follows we re-state the known examples where the conjugacy exists and, subsequently, give evidence of a more general pattern.

Let us start with $N=2$. The attractor of the standard IFS $f_j(z)=\lambda z+(-1)^j$ for $j=0,1$ is not invertible since $f_0^{-1}(z)\neq f_1^{-1}(z)$ for every $z\in\mathbb{C}$. However, the attractor is symmetric with respect to the origin, so $-A_\lambda=A_\lambda$. Equivalently, let $\sigma:\Sigma_2\to\Sigma_2$ be $\sigma(x)=x\bmod2$ and define $s_{\sigma(j)}(z)=(-1)^{\sigma(j)}z$, then $s_{\sigma(j)}(A_\lambda)=A_\lambda$ for $j=0,1$. From Lemma \ref{lem:rotationalsymmetry}, $A_\lambda$ is also the attractor of the IFS $\Phi_{2,\sigma}=\{f_j\circ s_{\sigma(j)}\}_{j=0}^1$, which implies that $q_\lambda$ is well defined for all $\lambda\in\mathcal{T}_2$, since $(f_0\circ s_{\sigma(0)})^{-1}(z)=(f_1\circ s_{\sigma(1)})^{-1}(z)$ exactly at $z=0$. For such $\lambda$ the attractor is a dendrite, which is connected and simply connected, and, if the address of $q_\lambda(0)$ is eventually periodic, Theorem \ref{thm:kameyamaPolynomial} applies and there is a homeomorphism between the dynamical systems $(A_\lambda,q_\lambda^n)$, for some integer $n$, and $(J(g),g)$ for some polynomial $g$, which has to be quadratic and can be chosen to be $Q_c(z)=z^2+c$. The integer $n$ has to be equal to $1$ since the action of $q_\lambda$ on $A_j$ is linear.
\begin{thm}[\cite{ErogluRohdeSolomyak2010}]
    \label{thm:boundedTurningQuadratic}
    Let $\lambda\in\mathcal{T}_2$ and suppose that the attractor $A_\lambda$ of the IFS $\Phi_{2,\sigma}$ is i.p.c.f. with $\sigma(x)=x\bmod2$. Then $A_\lambda$ has the bounded turning property.    
\end{thm}
Therefore, since dendritic Julia sets of $Q_c$ do not have parabolic cycles, the topological conjugacy between $(A_\lambda,q_\lambda)$ and $(J_c,Q_c)$ is actually quasisymmetric. The authors in \cite{ErogluRohdeSolomyak2010} revisits this conjugacy also form the symbolic dynamics point of view. They show
\begin{thm}
    \label{thm:SameKneadingSequence}
    Let $\lambda\in\mathcal{T}_2$ and $c\in\mathbb{C}$ is such that $J_c$ is a dendrite, and the kneading sequences of $q_\lambda$ and $Q_c$ are identical. Then $(A_\lambda,q_\lambda)$ and $(J_c,Q_c)$ are conjugate.
\end{thm}
Remember that the kneading sequence of $q_\lambda$ is defined as the itinerary of the critical value $q_\lambda(0)$ with respect to the partition $\{\{0\},A_1\setminus\{0\},A_0\setminus\{0\}\}$. Ero\u{g}lu \emph{et al.} describe how to obtain the kneading sequence from the address of $0$ under the maps $\{f_j\}_{j=0}^1$ and viceversa. We remark that the kneading sequence can be read directly from the neighbor graph $NG_\lambda$ of $\Phi_{2,\sigma}=\{f_0\circ s_{\sigma(0)},f_1\circ s_{\sigma(1)}\}$. 

\begin{ex}
Consider the attractor in the first row of Figure \ref{fig:attractorvsJulia}. It is taken from Example 5.4.1 in \cite{ErogluRohdeSolomyak2010}. The parameter $\lambda\approx 0.3668760+0.5202594i$ is the root of the power series \[1+\frac{-z+z^2+z^3}{1+z^3}=\frac{1-z+z^2+2z^3}{1+z^3}=\sum_{k=0}^\infty a_kz^k.\] The address of $0$ can be associated to the sequence of coefficients $\mathbf{a}=a_0a_1\cdots=+(-+++--)$, where the symbols inside the parentheses are repeated indefinitely. The kneading sequence $\nu=v_1v_2\cdots$ is obtained by reading consecutive pairs of elements from the sequence $\mathbf{a}$ and assign $v_j=\mathbf{0}$ if $a_{j-1}=a_j$ and $v_j=\mathbf{1}$ otherwise. In this example, we have $\nu=\mathbf{1}\overline{\mathbf{100}}$. Note that $\nu$ is an admissible kneading sequence since $\nu>\tau^k\nu$ for all $k\geq1$ where $\tau$ is the right-shift map and $>$ is the lexicographical ordering on $\Sigma_2^\infty$. To $\nu$ we can then associated an angle $\theta\in\mathbb{R}/2\pi\mathbb{Z}$ by simply flipping the symbols $\mathbf{0}$ and $\mathbf{1}$ and interpret it as a binary expansion: $\theta=(0.\mathbf{0}\overline{\mathbf{011}})_2=(1/8+1/16)(8/7)=3/14$. Since the angle is of the form $t2^{-k}(2^n-1)^{-1}$, the landing point of $R_{\mathcal{M}}(3/14)$ is a Misiurewicz parameter $c$, which is approximately $-0.155788+1.11222i$.

The neighbor graph of $\Phi_2=\{f_0,f_1\}$ at $\lambda$ is the first from the left in Figure \ref{fig:finitelyramified}. The $NG_\lambda$ for the invertible IFS $\Phi_{2,\sigma}$ is the leftmost in Figure \ref{fig:ipcfNeighborGraphs}. By reading the left symbol of the labels of the only infinite path on $NG_\lambda$, we recover the new addresses of $0$: $\mathbf{01}\overline{\mathbf{100}}$ and $\mathbf{11}\overline{\mathbf{100}}$. Consequently, the kneading sequence of the critical value under $q_\lambda$ is obtained by applying the right-shift map on either address. 

\begin{figure}
    \centering
    \begin{tikzpicture}
    
    \node [draw, circle] (0) at (0:0) (id) {\small $id$};
    \node [draw, fill=black, circle] (120) at (120:3) (s0a) {};
    \node [draw, fill=black, circle] (0) at (0:3) (s0c) {};
    \node [draw, fill=black, circle] (-60) at (-60:3) (s1a) {};
    \node [draw, fill=black, circle] (-120) at (-120:3) (s1b) {};

    \path[-{Triangle[]}] (id) edge [bend right, right] node {$01$} (s0a);
    \path[-{Triangle[]}] (id) edge [bend left, left] node {$10$} (s0a);
    
    \path[-{Triangle[]}] (s0a) edge [bend right, left] node {$11$} (s1b);
    \path[-{Triangle[]}] (s1b) edge [above] node {$11$} (s0c);
    \path[-{Triangle[]}] (s0c) edge [right] node {$00$} (s1a);
    \path[-{Triangle[]}] (s1a) edge [below] node {$00$} (s1b);
        
\end{tikzpicture}%
    \begin{tikzpicture}
    
    \node [draw, circle] (0) at (0:0) (id) {\small $id$};
    \node [draw, fill=black, circle] (-10) at (-10:2) (s0a) {};
    \node [draw, fill=black, circle] (10) at (10:2) (s0b) {};
    \node [draw, fill=black, circle] (100) at (110:3) (s1a) {};
    \node [draw, fill=black, circle] (135) at (135:2.5) (s1b) {};
    \node [draw, fill=black, circle] (250) at (250:2.5) (s2a) {};
    \node [draw, fill=black, circle] (210) at (210:4) (s2b) {};

    \path[-{Triangle[]}] (id) edge [bend left=10, above] node {$12$} (s0b);
    \path[-{Triangle[]}] (id) edge [bend right=10, below] node {$21$} (s0a);
    \path[-{Triangle[]}] (id) edge [bend right=0, right] node {$10$} (s1a);
    \path[-{Triangle[]}] (id) edge [bend left=0, left] node {$01$} (s1b);
    \path[-{Triangle[]}] (id) edge [bend left=10, right] node {$20$} (s2a);
    \path[-{Triangle[]}] (id) edge [bend right=0, right] node {$02$} (s2b);
    \draw[-{Triangle[]}] (s0a) to[out=-45, in=45, looseness=10] node[right] {$00$} (s0a);
    \draw[-{Triangle[]}] (s0b) to[out=-45, in=45, looseness=10] node[right] {$00$} (s0b);
    \path[-{Triangle[]}] (s1a) edge [bend right=100, left] node {$11$} (s2b);
    \path[-{Triangle[]}] (s2b) edge [bend left=150, left] node {$22$} (s1a);
    \path[-{Triangle[]}] (s1b) edge [bend right=100, left] node {$11$} (s2a);
    \path[-{Triangle[]}] (s2a) edge [bend left=150, left] node {$22$} (s1b);
    
\end{tikzpicture}
    \begin{tikzpicture}
    
    \node [draw, circle] (0) at (0:0) (id) {\small $id$};
    \node [draw, fill=black, circle] (00) at (0:3) (s0) {};
    \node [draw, fill=black, circle] (120) at (120:3) (s1) {};
    \node [draw, fill=black, circle] (240) at (240:3) (s2) {};

    \path[-{Triangle[]}] (id) edge [bend right=40, below] node {$10$} (s0);
    \path[-{Triangle[]}] (id) edge [bend right=10, below] node {$34$} (s0);
    \path[-{Triangle[]}] (id) edge [bend left=10, above] node {$43$} (s0);
    \path[-{Triangle[]}] (id) edge [bend left=40, above]node {$01$} (s0);
    \path[-{Triangle[]}] (id) edge [bend left=40, left] node {$20$} (s1);
    \path[-{Triangle[]}] (id) edge [bend left=10, left] node {$50$} (s1);
    \path[-{Triangle[]}] (id) edge [bend right=10, right] node {$05$} (s1);
    \path[-{Triangle[]}] (id) edge [bend right=40, right] node {$02$} (s1);
    \path[-{Triangle[]}] (id) edge [bend left=40, right] node {$12$} (s2);
    \path[-{Triangle[]}] (id) edge [bend left=10, right] node {$45$} (s2);
    \path[-{Triangle[]}] (id) edge [bend right=10, left] node {$54$} (s2);
    \path[-{Triangle[]}] (id) edge [bend right=40, left] node {$21$} (s2);
    
    \draw[-{Triangle[]}] (s0) edge [bend left=40, below] node {$00$} (s2);
    \draw[-{Triangle[]}] (s1) edge [bend left=40, right] node {$44$} (s0);
    \draw[-{Triangle[]}] (s2) edge [bend left=40, left] node {$22$} (s1);
    
\end{tikzpicture}
    \caption{Starting from the left, the neighbor graph for the dendrite of Figure \ref{fig:attractorvsJulia}, the Sierpi\'{n}ski triangle, and the Hexagasket when the IFS defining them is invertible. Compare these with the ones in Figure \ref{fig:finitelyramified}.}
    \label{fig:ipcfNeighborGraphs}
\end{figure}
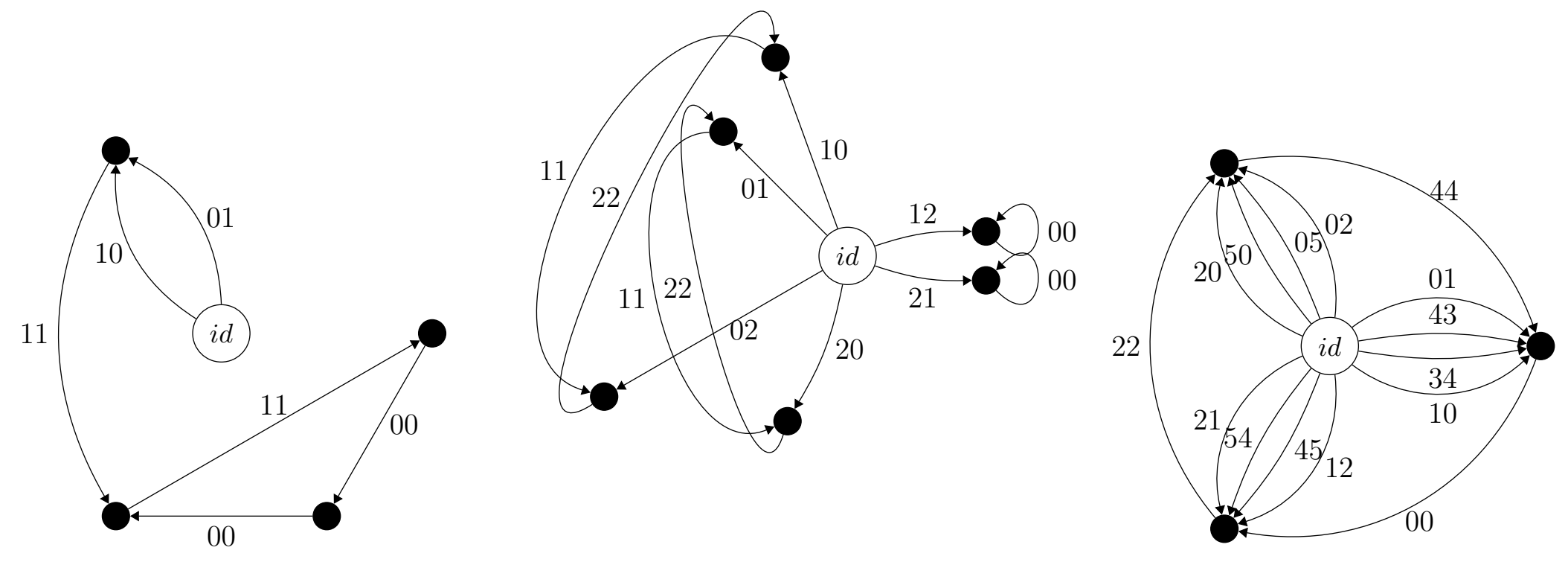
\end{ex}

When $N>2$, Theorem \ref{thm:SameKneadingSequence} can be extended to rational maps. Recall form Section \ref{sect:IFS} that $\lambda\in\mathbb{B}_N$ if and only if $\lambda$ is a root of a power series with coefficient in the set $\Delta_N=\{(\xi_N^a-\xi_N^b)/(1-\xi_N)~\mid~a,b\in\Sigma_N\}$. Let $\mathcal{P}_{\Delta_N}$ be the set of power series with coefficients in $\Delta_N$, normalized so that the constant term is $1$.
\begin{thm}[\cite{ErogluRohdeSolomyak2010}]
    \label{thm:boundedTurningHigherN}
    Assume the attractor $A_\lambda$ of the IFS $\Phi_N=\{f_j\}_{j=0}^{N-1}$ is p.c.f. In addition, suppose that for all $F\in \mathcal{P}_{\Delta_N}$ for which $F(\lambda)=0$ there exists $C_1>0$ such that all sufficiently large positive integers $n$ and $G\in\mathcal{P}_{\Delta_N}$ \[
    \lvert G(\lambda)\rvert\geq C_1\lvert\lambda\rvert^n\] whenever $n$ is the maximal common initial part of of $F$ and $G$. Then $A_\lambda$ has bounded turning.   
\end{thm}
If we add the assumption that the attractor is invertible, Theorem \ref{thm:kameyamaRational} applies and yield a p.c.f. rational map $g$. Such a map has no parabolic cycles and its critical points in $J(g)$ are stictly preperiodic; hence $g$ is semi-hyperbolic. Since the attractor of the IFS $\Phi_{N,\sigma}=\{f_j\circ s_{\sigma(j)}\}_{j=0}^{N-1}$ is i.p.c.f. with bounded turning, then the conjugacy is quasisymmetric by \cite[Theorem 3.1]{ErogluRohdeSolomyak2010}.

\begin{ex}
Set $N=3$. As before, any attractor obtained with the standard maps $f_j(z)=\lambda z+\xi_3^j$ for $j=0,1,2$ is not invertible. However, $\Phi_3(A_\lambda)$ has $2\pi/3$-rotational symmetry (Lemma \ref{lem:rotationalsymmetry}) and the translations $(f_0\circ s_{\sigma(0)})^{-1}$ and $(f_1\circ s_{\sigma(1})^{-1}$ agree on the point \[p=\frac{\xi_3^{\sigma(1)}-\xi_3^{\sigma(0)+1}}{\xi_3^{\sigma(1)}-\xi_3^{\sigma(0)}}\] other than $\infty$. Hence, the attractor of the IFS $\Phi_{3,\sigma}=\{f_j\circ s_{\sigma(j)}\}_{j=1}^2$ is invertible when $A_0\cap A_1=\{p\}$.

The intersection point of $A_0$ and $A_1$ in the Sierpi\'{n}ski triangle ($\lambda=1/2$) is easily found to be $(1+\xi_3)=-\xi_3^2$. Therefore, we need to solve
\[p=-\xi_3^2\iff \frac{\xi_3^{\sigma(1)}-\xi_3^{\sigma(0)+1}}{\xi_3^{\sigma(1)}-\xi_3^{\sigma(0)}}=-\xi_3^2\] which, if we assume $\sigma(x)=ax+b\bmod3$ for any $a,b\in\{0,1,2\}$, has solutions $\sigma(x)=2x+b\bmod3$ with $b=0,1,2$. In \cite{Kameyama2000} it is shown that for $b=0$ the i.p.c.f. Sierpi\'{n}ski triangle is homeomorphic to the Julia set of the semi-hyperbolic rational map $R_{2,1,-16/27}(z)=z^2-16/(27z)$. In \cite{ErogluRohdeSolomyak2010}, they show the dynamics of $q$ and $R_{2,1,-16/27}$ agree and, thus, concluding that the conjugacy is quasisymmetric. 

The second directed graph in Figure \ref{fig:ipcfNeighborGraphs} is the neighbor graph of $\Phi_{3,\sigma}(A_{0.5})$ where $\sigma(x)=2x\bmod3$. The combinatorial dynamics of $q_\lambda$ can be read off from the $NG_\lambda$ and coincides with that of the rational map: one critical point is prefixed with preperiod $1$; and the remaining two critical points enter the same cycle of period $2$ after one iteration, see also Figure \ref{fig:Sierpinki3dynamics}. The addresses of the critical values correspond to the kneading sequences: one critical value is fixed with kneading sequence $\overline{\mathbf{0}}$; while the other two critical values are mapped to each other and have kneading sequence $\overline{\mathbf{12}}$ and $\overline{\mathbf{21}}$.
\begin{figure}
    \centering
    \includegraphics[width=0.45\linewidth]{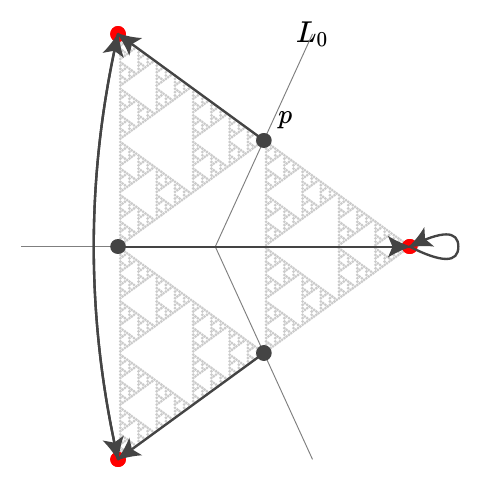}
    \includegraphics[width=0.4\linewidth]{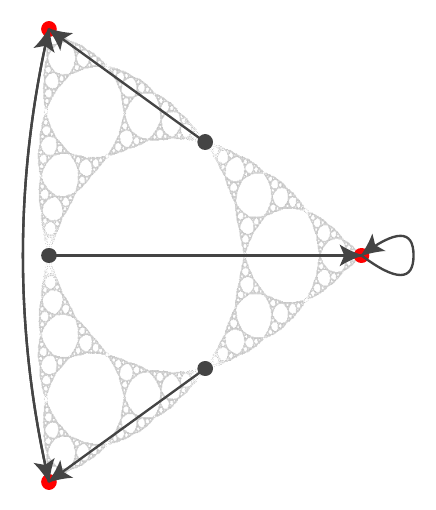}
    \caption{The dynamics of the critical points on the Sierpi\'{n}ski triangle and on the Julia set of the Misiurewicz rational map. The critical values are the corners of each set, while the critical points are the cut points.}
    \label{fig:Sierpinki3dynamics}
\end{figure}
Note that if we used the functions $\sigma(x)=2x+1\bmod3$ or $\sigma(x)=2x+2\bmod3$, then the diagram on the left in Figure \ref{fig:Sierpinki3dynamics} would be rotated by $2\pi/3$ and $4\pi/3$, respectively. Hence, the dynamics does not change, and the quasisymmetry conjugacy still holds. 
\end{ex}

\begin{ex}
The Hexagasket example was presented in \cite{ErogluRohdeSolomyak2010}. As the attractor of the standard IFS $\Phi_6=\{f_j\}_{j=0}^5$ with $\lambda=1/3$, the Hexagasket is not invertible. However, by Lemma \ref{lem:rotationalsymmetry} the attractor is symmetric wtih respect to a rotation by $2\pi/6$. The translations $(f_0\circ s_{\sigma(0)})^{-1}$ and $(f_1\circ s_{\sigma(1)})^{-1}$, for some function $\sigma$ on the set $\Sigma_6$, agree on the point \[p=\frac{\xi_6^{\sigma(1)}-\xi_6^{\sigma(0)+1}}{\xi_6^{\sigma(1)}-\xi_6^{\sigma(0)}}\] other than $\infty$. Simple calculations show that $A_0$ and $A_1$ intersect at the point $(1+\xi_6)/2$. Hence, assuming $\sigma(x)=ax+b\bmod6$ for $a,b=0,1,\cdots,5$, the solution to \[p=\frac{1+\xi_6}{2}\iff \frac{\xi_6^{\sigma(1)}-\xi_6^{\sigma(0)+1}}{\xi_6^{\sigma(1)}-\xi_3^{\sigma(0)}}=\frac{1+\xi_6}{2}\] is given by $\sigma(x)=3x+b\bmod6$. Setting $b=2$, the authors in \cite{ErogluRohdeSolomyak2010} showed that dynamics of $q_\lambda$ on the Hexagasket is equivalent to that of $R_{4,2,c}$ with $c=\tfrac{4}{9}\exp(2\pi i/6)$ on its Julia set. Considering that the Hexagasket has bounded turning and that the rational map is semi-hyperbolic, the two systems are quasisymmetrically conjugate. The two self-similar sets and the dynamics of the two maps are shown in Figure \ref{fig:Sierpinski6dynamics}. Recall that the Julia set of $R_{4,2,c}$ and $R_{4,2,\overline{c}}$ have conjugate dynamics \cite{DevaneyRochaSiegmund2007}, so the Julia set for $c=\tfrac{4}{9}\exp(2\pi i 5/6)$ is conjugate to that of $c=\tfrac{4}{9}\exp(2\pi i/6)$, and the outer arrows would move clockwise like the ones in the Hexagasket.
\begin{figure}
    \centering
    \includegraphics[width=0.45\linewidth]{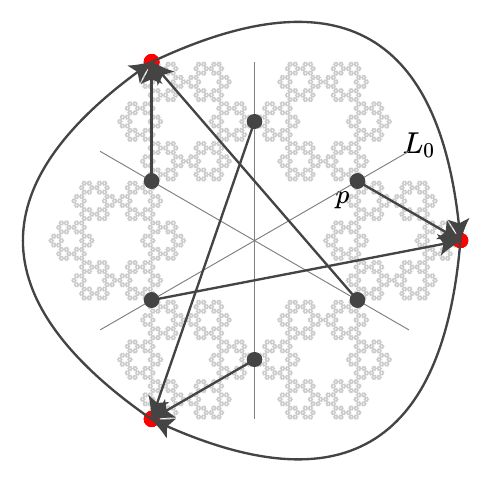}
    \includegraphics[width=0.45\linewidth]{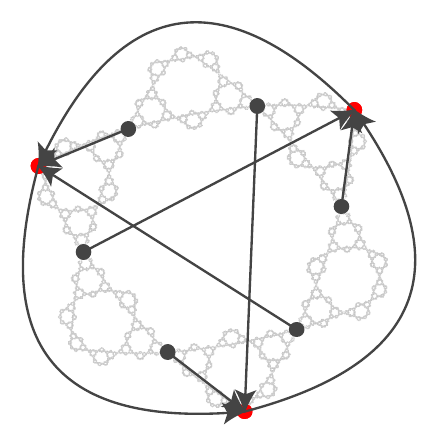}
    \caption{The dynamics of the critical points on the Hexagasket and on the Julia set of the Misiurewicz rational map. The critical values are the corners of each set, while the critical points are the cut points.}
    \label{fig:Sierpinski6dynamics}
\end{figure}

Observe also that if we use $\sigma(x)=3x+b\bmod6$ for $b=1,4,5$, we get conjugate orbit diagrams. If, instead, we set $b=0$ or $3$, the critical values now become three different fixed points. This different behavior of the critical points is consistent with the dynamics of $R_{n,m,c}$ with $c=-4/9$ whose Julia set is still a generalized Sierpi\'{n}ski gasket, see Figure \ref{fig:Sierpinski6dynamicsReal}.
\begin{figure}
    \centering
    \includegraphics[width=0.45\linewidth]{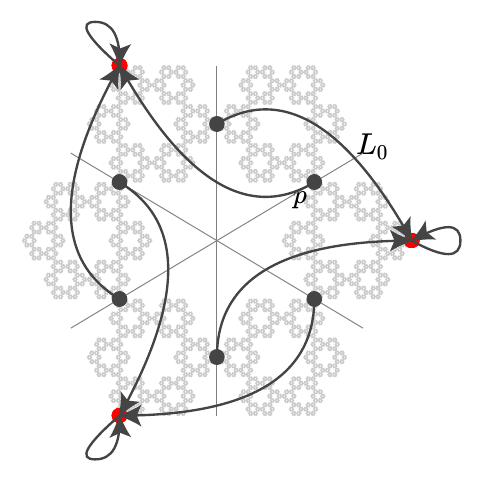}
    \includegraphics[width=0.45\linewidth]{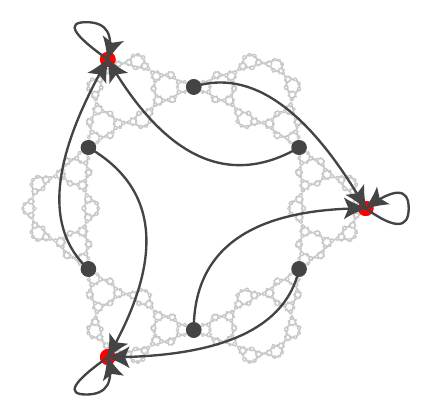}
    \caption{The dynamics of the critical points on the Hexagasket when $\sigma(x)=3x\bmod6$ and on the Julia set of $R_{4,2,-4/9}$.}
    \label{fig:Sierpinski6dynamicsReal}
\end{figure}
\end{ex}

We want to point out that the attractor $\Phi_6(A_\lambda)$ is always the Hexagasket when $\lambda=\tfrac{1}{3}\exp(2\pi ik/6)$ for $0\leq k\leq5$. Let $k=2$ and keep $\sigma(x)=3x\bmod6$, then we obtain the exact same orbit diagram of Figure \ref{fig:Sierpinski6dynamics}, or a conjugate one if we set $k=1$.

\begin{figure}
    \centering
    \includegraphics[width=0.45\linewidth]{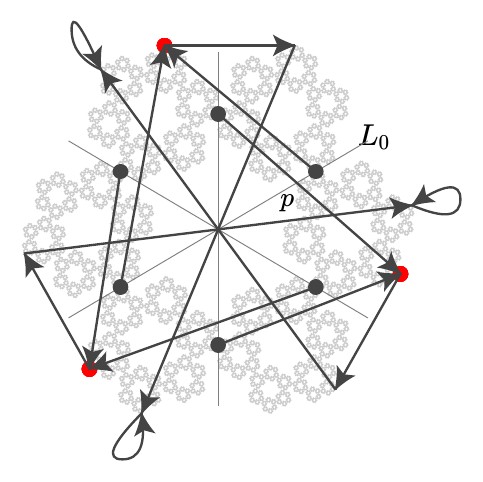}
    \includegraphics[width=0.45\linewidth]{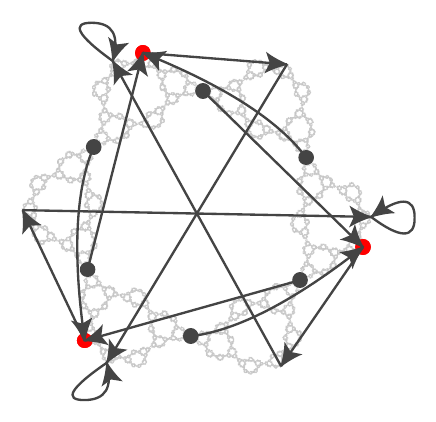}
    \caption{The dynamics of the critical points on $\Phi_6(A_\lambda)$ for $\lambda\approx0.338+0.08i$ when $\sigma(x)=3x\bmod6$ and on the Julia set of $R_{4,2,c}$ for $c\approx-0.4045+0.125i$.}
    \label{fig:Sierpinski6dynamicsNew}
\end{figure}

Essentially, we suggest to stick with one choice of $\sigma$, say $\sigma(x)=3x\bmod6$, then all possible connections between $\lambda\in\mathcal{T}_6$ and Misiurewicz parameters $c\in\mathcal{M}(R_{4,2,c})$ can be found. In fact, a quick computation gave us another parameter $\lambda\approx0.33791743+0.08029265i$ which is a root of \[\pi_\lambda(\mathbf{021}\overline{\mathbf{4}})=\pi_\lambda(\mathbf{121}\overline{\mathbf{4}})\implies(1-\xi_6)+2\xi_6^2\lambda+2\xi_6\lambda^2+\frac{2\xi_6\lambda^3}{1-\lambda}=0\]
and seems to be an element of $\mathcal{T}_6$, see Figure \ref{fig:Sierpinski6dynamicsNew}. The corresponding Misiureiwicz parameter in the McMullen family $R_{4,2,c}$, if any, should have the same kneading sequence: considering that the critical points are eventually periodic and the position of the critical values are roughly the same as those of $\lambda=1/3$, then an educated guess seems to be $c\approx-0.4045230912763623+0.12470417452497898i$.

We have not seen in the literature any other correspondence between parameters $\lambda$ and $c$, so we provide possibly new ones. For starters, when $N\not\equiv0\bmod4$ in \cite{BandtHung2008b} gave the minimum absolute value for $\lambda\in\mathbb{D}$ to guarantee that the attractor $\Phi_N(A_\lambda)$ is ``just touching'', i.e. that $\#(A_0\cap A_1)=1$:
\[\lvert \lambda\rvert=\begin{cases}
    \frac{\sin(\pi/N)}{1+\sin(\pi/N)} &N\equiv 2\bmod4\\
    \frac{\sin(\pi/N)}{\cos(\pi/(2N))+\sin(\pi/N)} &N\equiv 1,3\bmod4
\end{cases}\]
If $\lambda\in\mathbb{R}$ and equal to the value above, we obtain a generalized Sierpi\'{n}ski gasket. Furthermore, we claim we can always find a corresponding real Misiurewicz parameter $c$ whenever $N\not\equiv0\bmod4$. When $N=4$ we also seem to have found cases where the conjugacy seems to exists. The table below records the ``easiest'' pair of parameters $\lambda$ and $c$ that we have found or have been found.

\begin{table}[ht]
\caption{Summary of the proved and guessed correspondence between IFS parameters and Misurewicz parameters of McMullen maps. The proved ones are those when $N=2,3,6$. When $N=2$, the McMullen map $R_{2,0,c}$ is, in fact, $Q_c$.}
{\centering
\begin{tabular}{ccccc}
\hline
$N$ & $\sigma(x)$ & $(n,m)$ & $\lambda$ & $c$ \\
\hline
$2$ & $x$   & $(2,0)$ & $1/2$             & $-2$ \\
$3$ & $2x$  & $(2,1)$ & $1/2$             & $-16/27$ \\
$4$ & $2x$  & $(3,1)$ & $(2+i)/5$         & $\approx -0.1875+0.5303i$ \\
$5$ & $2x$  & $(4,1)$ & $\sin(\pi/5)/(\cos(\pi/10)+\sin(\pi/5))$           & $\approx -0.6307528$ \\
$6$ & $3x$  & $(4,2)$ & $1/3$             & $-4/9$ \\
$7$ & $4x$  & $(4,3)$ & $\sin(\pi/7)/(\cos(\pi/14)+\sin(\pi/7))$           & $\approx-0.3455958$ \\$8$ & ?  & ? & ?           & ? \\
$9$ & $4x$  & $(6,3)$ & $\sin(\pi/9)/(\cos(\pi/18)+\sin(\pi/9))$           & $\approx-0.41959445$ \\
$10$ & $5x$  & $(6,4)$ & $\sin(\pi/10)/(1+\sin(\pi/10))$           & $\approx-0.35336916$ \\
$11$ & $6x$  & $(6,5)$ & $\sin(\pi/11)/(\cos(\pi/22)+\sin(\pi/11))$           & $\approx-0.30567227$ \\
$12$ & ? & ? & ? & ? \\
$13$ & $6x$  & $(8,5)$ & $\sin(\pi/13)/(\cos(\pi/26)+\sin(\pi/13))$           & $\approx-0.359237706$ \\
\hline
\end{tabular}}
\label{tbl:correspondenceSummary}
\end{table}
Note that when $N=4$ the smallest real parameter $\lambda$ in $\mathbb{B}_N$ is $0.5$, but in that case the attractor is a square. In fact, when $N\equiv0\bmod4$ the smallest real value of $\lambda$ in $\mathbb{B}_N$ is given \cite{Strichartz1999} by the formula \[\frac{\sin(\pi/N)}{\cos(\pi/N)+\sin(\pi/N)}\] and the intersection $A_0\cap A_1$ is a Cantor set except for $N=4$, where, instead, it is a segment of $L_0$. We were unable to find an easy pair of parameters for $N=8$ but most probably the tuple $(n,m)$ should be $(5,3)$ and $\sigma(x)$ should be $4x\bmod8$. Similarly for $N=12$ we believe $(n,m)=(7,5)$ and $\sigma(x)=6x\bmod12$. In both cases, $c$ would not be real, just like when $N=4$.
\begin{quest}
    Is it possible to find a description of the subset of $\mathcal{T}_N$ for $N>2$ for which the attractor is i.p.c.f.?
\end{quest}
Thinking back to Lemma \ref{lem:ipcfsubsetofT} and to the examples in Section \ref{subsect:propertiesAttractor} of attractors which are not i.p.c.f., we prove the following:
\begin{prop}
    Let $N\geq2$ and suppose the attractor $A_\lambda$ of the IFS $\Phi_N$ is p.c.f. with $A_k\cap A_{k+1}=\{p_k\}$ (and no other intersections between the $A_j$) for each $0\leq k<N$. Assume further that there exists $0<M<N$ such that for each $k$ the pair of addresses $(\mathbf{k}a_{k,1}a_{k,2}\cdots,(\mathbf{k+1})b_{k,1}b_{k,2}\cdots)$ in $\pi_\lambda^{-1}(p_k)$ satisfy $a_{k,j}\equiv b_{k,j}-M\bmod N$ for all $j\geq1$. Then $A_\lambda$ is i.p.c.f. as the attractor of the IFS $\Phi_{N,\sigma}$ where $\sigma(x)=Mx\mod N$.
\end{prop}
\begin{proof}
    To show that $A_\lambda$ is i.p.c.f. is enough to prove that for each $0\leq k<N$
    \[(f_k\circ s_{\sigma(k)})^{-1}(p_k)=(f_{k+1}\circ s_{\sigma(k+1)})^{-1}(p_k)\iff f_{k}^{-1}(p_k)=\xi_N^{\sigma(k)-\sigma(k+1)}f_{k+1}^{-1}(p_k).\]
    Indeed, by assumption $a_{k,j}\equiv b_{k,j}-m\bmod N$ and $\sigma(k)-\sigma(k+1)=-M$, so 
    \begin{align*}
        f_{k}^{-1}(p_k)=\pi_\lambda(a_{k,1}a_{k,2}\cdots)&=\sum_{j=1}^{\infty}\xi_N^{a_{k,j}}\lambda^{j-1}=\sum_{j=1}^{\infty}\xi_N^{b_{k,j}-M}\lambda^{j-1}=\xi_N^{-M}\sum_{j=1}^{\infty}\xi_N^{b_{k,j}}\lambda^{j-1}\\
        &=\xi_N^{-M}\pi_\lambda(b_{k,1}b_{k,2}\cdots)=\xi_N^{-M}f_{k+1}^{-1}(p_k)
    \end{align*}
    as claimed. 
\end{proof}

\begin{quest}
Fix $N\geq2$ and consider the subset of $\mathcal{T}_N$ for which the attractor is i.p.c.f.. Is the pair $(n,m)$ constant for all parameters in this subset? If not, what are all the possible pairs $(n,m)$?
\end{quest}
For example, the Julia sets of the rational function $R_{1,1,c}$ are \cite{DevaneyHolzerUminsky2006} either Cantor sets, or the whole sphere, or are Jordan curves passing through $0$ and $\infty$. The Julia set of $R_{0,2,c}$ instead is always a circle when $c\neq0$. Hence, only the pair $(n,m)=(2,0)$ is of our interest and $\mathcal{T}_2$ is related to $\mathcal{M}(Q_c)$. Analogous reasons tell us that $\mathcal{T}_N$ is not related to $\mathcal{M}(R_{1,m,c})$ (see \cite{CampbellCollins2017}), nor to $\mathcal{M}(R_{0,m,c})$ for $m\geq3$.

\begin{figure}
    \centering
    \includegraphics[width=0.45\linewidth]{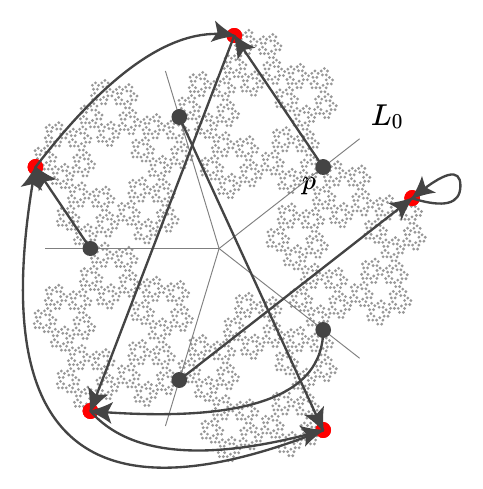}
    \includegraphics[width=0.45\linewidth]{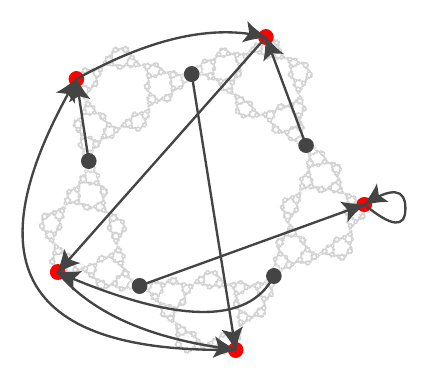}
    \caption{The dynamics of the critical points on $\Phi_5(A_\lambda)$ for $\lambda\approx0.37027+0.1525i$ when $\sigma(x)=3x\bmod5$ and on the Julia set of $R_{3,2,c}$ for $c\approx-0.134355+0.353877i$.}
    \label{fig:Sierpinski5dynamicsNew}
\end{figure}

Table \ref{tbl:correspondenceSummary} shows that for at least $N\leq13$, there seems to be a pattern in the choice of the pair $(n,m)$ relating $\mathcal{T}_N$ to $\mathcal{M}(R_{n,m,c})$.  
\begin{conj}
    Set integer $N=4k+r>2$ and let $\lambda\in\mathcal{T}_N\cap\mathbb{R}$. Suppose further that the attractor $A_\lambda$ of the IFS $\Phi_{N,\sigma}$ is i.p.c.f. for some $\sigma:\Sigma_N^\infty\to\Sigma_N^\infty$. Then there exists Misiurewicz parameter $c$ in $\mathcal{M}(R_{n,m,c})$ with $N=n+m$ such that the dynamical systems $(A_\lambda,q_\lambda)$ and $(J_c,R_{n,m,c})$ are quasisymmetrically conjugate. Specifically, we have that $\sigma(x)=(m+1)x\bmod N$ and if $N\equiv0\bmod 4$ then $n=2k+1$ and $m=2k-1$, otherwise $n=2k+2$ and $m=2k-2+r$. Moreover, $c\in\mathbb{R}\cap\mathcal{M}(R_{n,m,c})$ if $N\not\equiv0\bmod4$.
\end{conj}
Is Lemma \ref{lem:mcorners}, which determines precisely the number of corners that lie on an edge of the boundary of the trap door in $J(R_{n,m,c})$ at each iteration, the reason behind the formulation of $\sigma$ in the above conjecture?
\begin{quest}
    What is the main reason behind the failing of some Misiurewicz parameters $c\in\mathcal{M}(R_{n,m,c})$ or $c\in\mathcal{M}(Q_c)$ to have a corresponding $\lambda\in\mathcal{T}_N$?
\end{quest}
Can we find other pairs of $(n,m)$ for larger $N$? Consider, for example, $(n,m)=(3,2)$ at the Misiurewicz parameter $c\approx -0.13435514872957918+0.35387748202031943i$. The kneading sequence of a critical point is $\overline{\mathbf{1342}}$. Thus, imposing the equation $\pi_\lambda(\mathbf{0}\overline{\mathbf{1342}})=\pi_\lambda(\mathbf{1}\overline{\mathbf{1342}})$ with $\sigma(x)=(m+1)x\bmod 5$, the resulting parameter is $\lambda\approx0.37026809+0.152509i$. A plot of the orbit diagram of the critical points of $q$ and that of $R_{3,2,c}$ is shown in Figure \ref{fig:Sierpinski5dynamicsNew}. In order to be thorough, we also tried different $\sigma$ when solving the equation but, for all of the resulting $\lambda$, the attractor had a significant overlap set. 

Following the same reasoning, we considered $(n,m)=(5,2)$ at the Misiurewicz parameter $c\approx-0.30243818206484524+0.372130731889976i$; found the kneading sequence of a critical point to be $\overline{\mathbf{2}}$; solved the equation $\pi_\lambda(\mathbf{0}\overline{\mathbf{2}})=\pi_\lambda(\mathbf{1}\overline{\mathbf{2}})$ with $\sigma(x)=(m+1)x\bmod 7$ and obtained the parameter $\lambda\approx0.26222933+0.17397387i$. 

\section{Conclusion}
\label{sect:Conclusion}
In this last section we indicate the tools that in our view might be useful in answering the raised questions and conjectures of the previous sections.

\subsection{Hubbard trees, geodesics and sprouts}
The case when $N=2$ is where the most is already in place, both in the IFS setting and in the corresponding quadratic case $Q_c(z)=z^2+c$. Hubbard trees have been studied extensively, both combinatorially and metrically, and Ero\u{g}lu, Rohde, and Solomyak already show that some of their characteristics are mirrored by the corresponding tree structure on the IFS side, see Section \ref{subsect:raysKSentropy}. 

What is the Hubbard tree in the IFS setting? By definition, a Hubbard tree is the minimal invariant tree connecting the critical point to elements in its orbit. Such object is constructible only in an i.p.c.f fractal set via what Bandt \cite{Bandt2014,BandtStahnke1990,Strichartz1999} calls \emph{geodesics}: minimal paths connecting (dynamical) boundary points in the attractor. These are constructed via a graph directed construction, much like the neighbor graph. His idea is that once these finitely many geodesics $\mathcal{C}_j$ are defined, a countable union of copies of them can be used to define arcs connecting arbitrary points $x,y$ in the attractor and, in this way, compute the interior distance $\rho(x,y)$.

\begin{figure}
    \centering
    \includegraphics[width=0.325\linewidth]{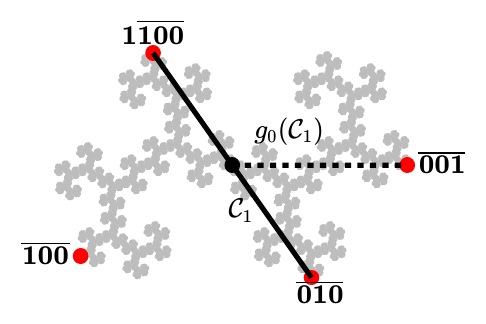}%
    \includegraphics[width=0.325\linewidth]{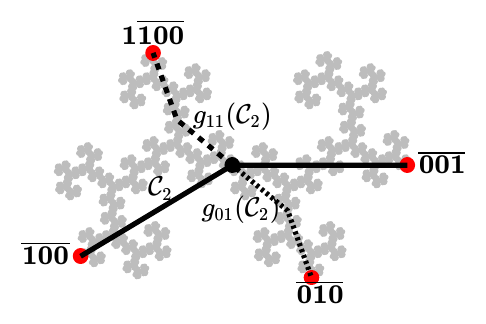}%
    \includegraphics[width=0.325\linewidth]{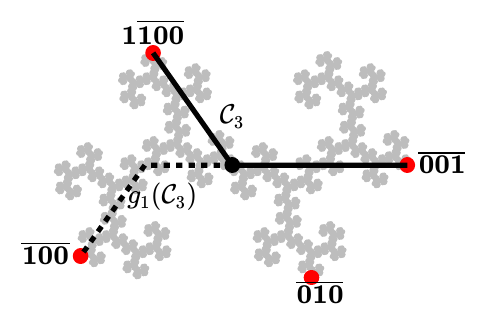}
    
    \includegraphics[width=0.325\linewidth]{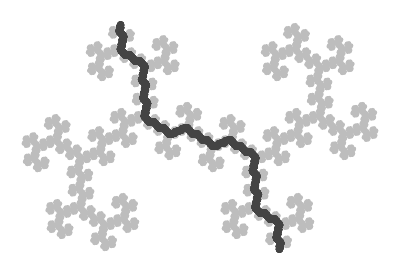}%
    \includegraphics[width=0.325\linewidth]{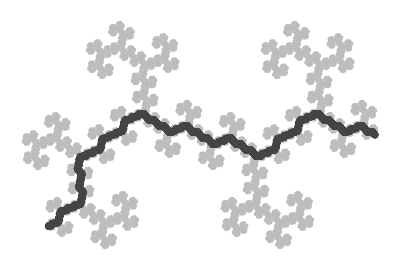}%
    \includegraphics[width=0.325\linewidth]{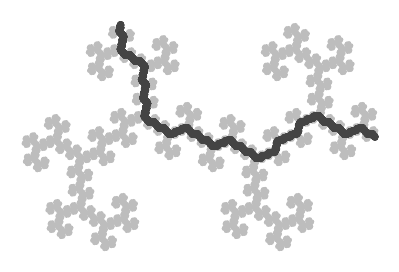}
    
    \caption{In the figure of the top row, the solid lines represent first level approximations of the geodesics. Their images are indicated with dashed or dotted lines. The bottom row are the actual geodesics.}
    \label{fig:geodesics}
\end{figure}

\begin{figure}[!htb]
    \centering
    \includegraphics[width=0.45\linewidth]{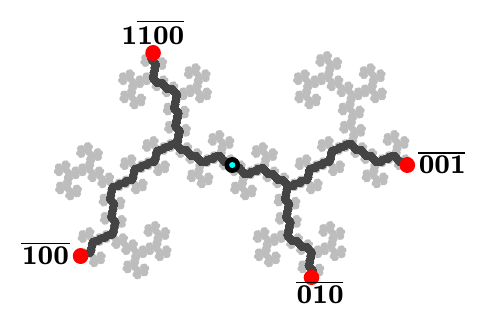}%
    \includegraphics[width=0.42\linewidth]{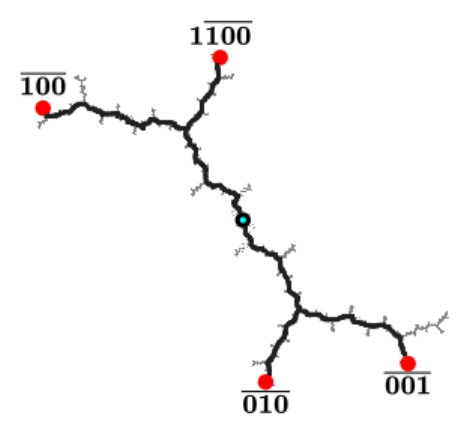}
    \caption{The union of the three geodesics can be thought of as the counterpart of the Hubbard tree in the IFS setting. On the right, the Hubbard tree at the corresponding Misiurewicz parameter.}
    \label{fig:ifs2HubbardTree}
\end{figure}

We illustrate the method to obtain these geodesics via the first example in Section \ref{subsubsect:knownexamples}: set $N=2$, then the attractor $\Phi_2(A_\lambda)$ for the parameter $\lambda\approx0.366+0.520i$ is a dendrite on which the expanding map $q_\lambda$, defined as the inverse branches of $g_j=f_j\circ s_{\sigma(j)}$ for $j=0,1$ and $\sigma(x)=x\bmod N$, is well defined. Let $[z_1,z_2]\subset X$ denote a connected set of points $z_1,z_2\in X$. We want to consider curves connecting pair of points in the post critical set, i.e. the dynamical boundary of $A_\lambda$: 
\begin{align*}
\mathcal{C}_1&:=\{w\in A_\lambda~\mid~w\in[\pi_\lambda(\mathbf{1}\overline{\mathbf{100}}),\pi_\lambda(\overline{\mathbf{010}})]\},\\
\mathcal{C}_2&:=\{w\in A_\lambda~\mid~w\in[\pi_\lambda(\overline{\mathbf{100}}),\pi_\lambda(\overline{\mathbf{001}})]\},\\ 
\mathcal{C}_3&:=\{w\in A_\lambda~\mid~w\in[\pi_\lambda(\mathbf{1}\overline{\mathbf{100}}),\pi_\lambda(\overline{\mathbf{001}})]\}
\end{align*}
Notice how each of these ``arc'' has to go through the critical point $0$ which has the pair of addresses $(\mathbf{01}\overline{\mathbf{100}},\mathbf{11}\overline{\mathbf{100}})$. Then one can check, by following the endpoints or visually by the approximations in the top row of Figure \ref{fig:geodesics}, that 
\begin{align*}
    \mathcal{C}_1&=g_{11}(\mathcal{C}_2)\cup g_{01}(\mathcal{C}_2)\\
    \mathcal{C}_2&=g_{0}(\mathcal{C}_1)\cup g_{1}(\mathcal{C}_3)\\
    \mathcal{C}_3&=g_{0}(\mathcal{C}_1)\cup g_{11}(\mathcal{C}_2)
\end{align*}
The geodesics are shown in the bottom row Figure \ref{fig:geodesics}. Their union can be thought of as the ``Hubbard tree'' for this IFS. Compare in Figure \ref{fig:ifs2HubbardTree} this ``linear'' Hubbard tree with the ``quadratic'' one of $z^2+c$, where $c$ is the landing point of the external ray $R_{\mathcal{M}}(3/14)$.

\begin{figure}
    \centering
    \includegraphics[width=0.425\linewidth]{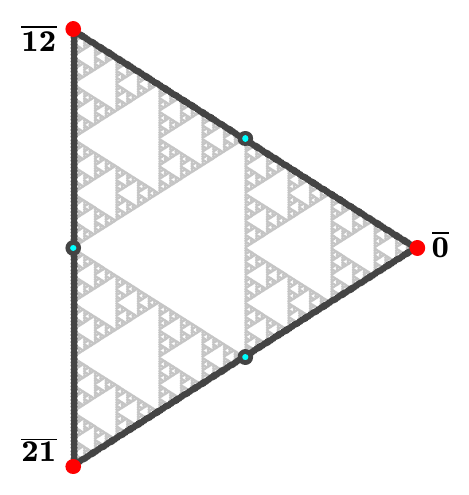}
    \includegraphics[width=0.425\linewidth]{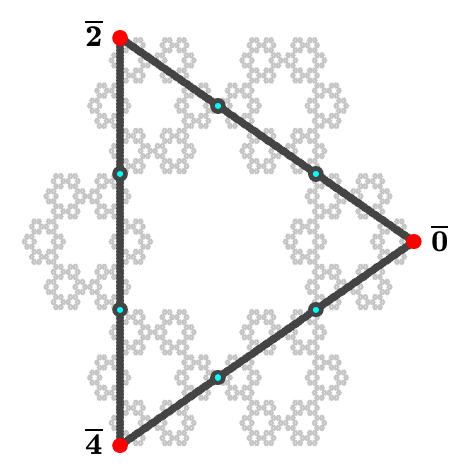}
    \caption{The geodesics for the Sierpi\'{n}ski triangle and Hexagasket are sides of an equilateral triangle. Compare with Figure \ref{fig:Sierpinki3dynamics} and Figure \ref{fig:Sierpinski6dynamicsReal}.}
    \label{fig:geodesics3and6}
\end{figure}

Similar reasoning can be made for the Siepri\'{n}ski triangle $(N,\lambda,\sigma(x))=(3,1/2,2x)$ and the Hexagasket $(N,\lambda,\sigma(x))=(6,1/3,3x)$. In both cases the geodesics turn out to be edges of an equilateral triangle, see Figure \ref{fig:geodesics3and6}. Specifically, for $N=3$ the geodesics are the unique compact sets $\mathcal{C}_j$ satisfying 
\begin{align*}
    \mathcal{C}_1&=g_{0}(\mathcal{C}_1)\cup g_{1}(\mathcal{C}_2)\\
    \mathcal{C}_2&=g_{2}(\mathcal{C}_1)\cup g_{1}(\mathcal{C}_3)\\
    \mathcal{C}_3&=g_{2}(\mathcal{C}_2)\cup g_{0}(\mathcal{C}_3)
\end{align*}
where $g_j:=f_{j}\circ s_{\sigma(j)}$ with $0\leq j\leq 2$ and $\sigma(x)=2x\bmod 3$; while for $N=6$ the geodesics satisfy
\begin{align*}
    \mathcal{C}_1&=g_{0}(\mathcal{C}_1)\cup g_{1}(\mathcal{C}_1)\cup g_{2}(\mathcal{C}_1)\\
    \mathcal{C}_2&=g_{2}(\mathcal{C}_2)\cup g_{3}(\mathcal{C}_2)\cup g_{4}(\mathcal{C}_2)\\
    \mathcal{C}_3&=g_{4}(\mathcal{C}_3)\cup g_{5}(\mathcal{C}_3)\cup g_{0}(\mathcal{C}_3)
\end{align*}
where $g_j:=f_{j}\circ s_{\sigma(j)}$ with $0\leq j\leq 5$ and $\sigma(x)=3x\bmod 6$.

\begin{figure}
    \centering
    \includegraphics[width=0.4\linewidth]{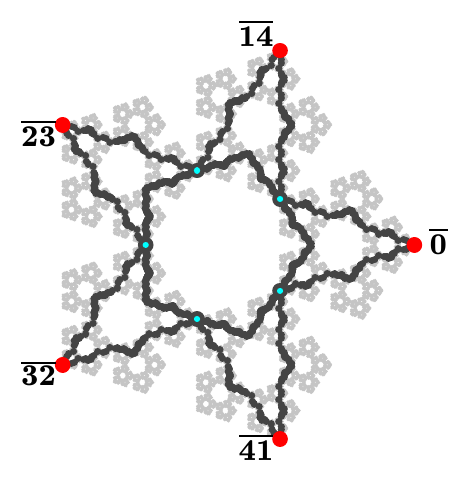}%
    \includegraphics[width=0.4\linewidth]{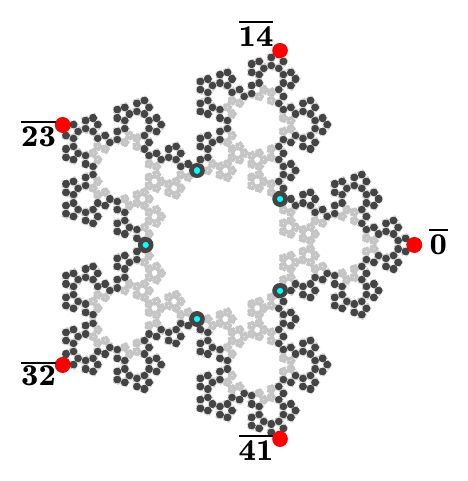}
    \caption{The union of five connecting curves between the points in the postcritical set of $q_\lambda$ when $N=5$. The connecting dimension of the union of the curves on the left is lower than the one the right.}
    \label{fig:geodesicsPentagasket}
\end{figure}

As pointed out in \cite{Strichartz1999}, there is a problem of ambiguity when defining the connecting curves $\mathcal{C}_j$. A priori there may be many ways to connect the dynamical boundary points through adjacent critical points and we do not know which one or which ones correspond to the minimal path. In his article, Strichartz is concerned with the \emph{connectivity dimension} $\delta$, i.e. the minimal $\delta$-dimensional Hausdorff measure of a connecting set, so he resolves the ambiguity by computing $\delta$ for all possibilities and choosing the set of $\mathcal{C}_j$ that give the smallest value. 
\begin{quest}
    Does the quasisymmetric image of the Hubbard tree in $J_c$ for Misiurewicz parameter $c$ correspond to the union of connecting sets $\mathcal{C}_j$ with minimal connectivity dimension in $\Phi_{2,\sigma}(A_\lambda)$ for $\lambda\in\mathcal{T}_2$?
\end{quest}
An example of such ambiguity is shown in Figure \ref{fig:geodesicsPentagasket} where $(N,\lambda,\sigma(x))$ with $N=5$ is chosen from Table \ref{tbl:correspondenceSummary}. The connecting curves on the left picture satisfy the following set of equations:
\begin{align*}
    \mathcal{C}_1=\{w\in \Phi_5(A_\lambda)~\mid~w\in[\pi_\lambda(\overline{\mathbf{0}}),\pi_\lambda(\overline{\mathbf{23}})]\}&=g_0(\mathcal{C}_1)\cup g_{12}(\mathcal{C}_2)\cup g_{11}(\mathcal{C}_5)\cup g_2(\mathcal{C}_2)\\
    \mathcal{C}_2=\{w\in \Phi_5(A_\lambda)~\mid~w\in[\pi_\lambda(\overline{\mathbf{14}}),\pi_\lambda(\overline{\mathbf{32}})]\}&=g_3(\mathcal{C}_1)\cup g_{20}(\mathcal{C}_1)\cup g_{21}(\mathcal{C}_3)\cup g_1(\mathcal{C}_5)\\
    \mathcal{C}_3=\{w\in \Phi_5(A_\lambda)~\mid~w\in[\pi_\lambda(\overline{\mathbf{23}}),\pi_\lambda(\overline{\mathbf{41}})]\}&=g_2(\mathcal{C}_4)\cup g_{34}(\mathcal{C}_2)\cup g_{30}(\mathcal{C}_4)\cup g_4(\mathcal{C}_5)\\
    \mathcal{C}_4=\{w\in \Phi_5(A_\lambda)~\mid~w\in[\pi_\lambda(\overline{\mathbf{32}}),\pi_\lambda(\overline{\mathbf{0}})]\}&=g_3(\mathcal{C}_3)\cup g_{43}(\mathcal{C}_3)\cup g_{44}(\mathcal{C}_5)\cup g_0(\mathcal{C}_4)\\
    \mathcal{C}_5=\{w\in \Phi_5(A_\lambda)~\mid~w\in[\pi_\lambda(\overline{\mathbf{41}}),\pi_\lambda(\overline{\mathbf{14}})]\}&=g_4(\mathcal{C}_2)\cup g_{03}(\mathcal{C}_1)\cup g_{02}(\mathcal{C}_4)\cup g_1(\mathcal{C}_3).
\end{align*}
Here each connecting curve goes through two critical points and connects boundary points whose ratio is $\xi_5^2$. Hence, two intersecting curves $\mathcal{C}_j$ do so either only at one point which is either one of their endpoints or is a critical points. The connecting curve on the right picture of Figure \ref{fig:geodesicsPentagasket} satisfy:
\begin{align*}
    \mathcal{C}_1=\{w\in \Phi_5(A_\lambda)~\mid~w\in[\pi_\lambda(\overline{\mathbf{0}}),\pi_\lambda(\overline{\mathbf{14}})]\}&=g_0(\mathcal{C}_1)\cup g_{0}(\mathcal{C}_2)\cup g_{1}(\mathcal{C}_3)\cup g_1(\mathcal{C}_4)\\
    \mathcal{C}_2=\{w\in \Phi_5(A_\lambda)~\mid~w\in[\pi_\lambda(\overline{\mathbf{14}}),\pi_\lambda(\overline{\mathbf{23}})]\}&=g_1(\mathcal{C}_1)\cup g_{2}(\mathcal{C}_2)\cup g_{2}(\mathcal{C}_3)\cup g_1(\mathcal{C}_5)\\
    \mathcal{C}_3=\{w\in \Phi_5(A_\lambda)~\mid~w\in[\pi_\lambda(\overline{\mathbf{23}}),\pi_\lambda(\overline{\mathbf{32}})]\}&=g_3(\mathcal{C}_1)\cup g_{3}(\mathcal{C}_2)\cup g_{2}(\mathcal{C}_4)\cup g_2(\mathcal{C}_5)\\
    \mathcal{C}_4=\{w\in \Phi_5(A_\lambda)~\mid~w\in[\pi_\lambda(\overline{\mathbf{32}}),\pi_\lambda(\overline{\mathbf{41}})]\}&=g_4(\mathcal{C}_1)\cup g_{3}(\mathcal{C}_3)\cup g_{3}(\mathcal{C}_4)\cup g_4(\mathcal{C}_5)\\
    \mathcal{C}_5=\{w\in \Phi_5(A_\lambda)~\mid~w\in[\pi_\lambda(\overline{\mathbf{41}}),\pi_\lambda(\overline{\mathbf{0}})]\}&=g_4(\mathcal{C}_2)\cup g_{4}(\mathcal{C}_3)\cup g_{0}(\mathcal{C}_4)\cup g_0(\mathcal{C}_5).
\end{align*}
Each $\mathcal{C}_j$ in this case contains only one critical point and connects (rotationally) adjacent boundary points. Two intersecting curve $\mathcal{C}_j$ do so only at their endpoints.
One can check that the connecting dimension of the connecting curves on the left is lower than the one on the right. 
\begin{quest}
    Consider the dynamical systems $(\Phi_{N,\sigma}(A_\lambda),q_\lambda)$ and $(J_c,R_{n,m,c})$ with parameters taken from Table \ref{tbl:correspondenceSummary}. Can the quasisymmetric image of the union of geodesics in $\Phi_{N,\sigma}(A_\lambda)$ be considered the Hubbard graph of in the Julia set of $R_{n,m,c}$ ?
\end{quest}

Tetenov, Yudin, and Drozdov \cite{TetenovYudinDrozdov2026} have constructed combinatorial schemes called \emph{sprouts} for any self-similar sets satisfying the finite intersection property, that is those satisfying $\#(A_j\cap A_i)<\infty$. These schemes carry the topological data of the attractor very much like the neighbor graph do. Each sprout $\Gamma$ has associated an \emph{index diagram}, a directed graph on the dynamical boundary of the attractor $A_\lambda$, which is built to record how a fixed correspondence between pieces and boundary acts. This is precisely what the induced map $q_\lambda$ supplies once $A_\lambda$ is i.p.c.f.: each directed edge of the index diagram, labeled by the map $f_j$ sending one boundary point to another, records one step of $q_\lambda$ acting on $\partial A_\lambda$. 
\begin{quest}
    For $\lambda\in\mathcal{T}_N$ such that $\Phi_{N,\sigma}(A_\lambda)$ is i.p.c.f. with respect to some $\sigma$, is the neighbor graph of $\Phi_{N,\sigma}$ isomorphic to the sprout of the same system? 
\end{quest}
The importance of such question comes from the fact that Tetenov, Yudin, and Drozdov show that these schemes determine self-similar dendrites up to isomorphism. Moreover, they provide an algorithm that enables the reconstruction of the main topological elements of a self-similar dendrite from its sprout. The algorithm uses what they call the \emph{main tree}, that is the minimal dendrite contained in the attractor such that it contains the dynamical boundary of the attractor. In other words, the main tree is the union of the main geodesics of Bandt when $N=2$.

\subsection{Core entropy}
Another reason to understand the counterpart of an Hubbard tree in the IFS setting is given by the core entropy. The monotonicity of core entropy along veins, proved in the post-critically finite setting by Li \cite{Li2007} and in general by subsequent work of Zeng \cite{Zeng2020}, is used to relate the order of parameters along a principal vein to the corresponding order of their growth rates. This monotonicity enters explicitly in the proof of Lemma 9.6 in \cite{LindseyTiozzoWu2025} and in bounding the range of growth rates on a principal vein.
\begin{quest}
    Let $R_\mathcal{M}(\theta)$ land at a Misiurewicz parameter $c\in\mathcal{M}(Q_c)$. Suppose that there is no corresponding $\lambda\in\mathcal{T}_2$. Order external angles along a vein according to the order of their landing points. Is it true that the absence of a correspondence persists for all $\theta'\prec \theta$ along the vein in $\mathcal{M}(Q_c)$ containing $c$?
\end{quest}
A positive answer would give a structural explanation for the experimental pattern observed in Figure \ref{fig:entropy}.

\subsection{Laminations.} 
Calegari and Walker \cite{CalegariWalker2026} develop a theory of laminations for the IFS $\Phi_2$. For parameters in $\lambda\in\partial\mathbb{B}_2$ they define the \emph{full set} $K$ of $A_\lambda$ to be the smallest compact simply-connected subset of $\mathbb{C}$ containing $A_\lambda$. They discover that the topology of $K$ is entirely captured by the structure of its set of cut points, and give a sufficient condition for the existence of cut points in $K$ for $\lambda\in\partial\mathbb{B}_2$. In fact, such condition is also necessary provided that a particular conjecture holds. The authors can associate to any parameter $\lambda\in\partial\mathbb{B}_2$ a tuple $(\theta,s)$ consisting of \emph{external angle} $\theta\in\mathbb{R}/2\pi\mathbb{Z}$ and a scale factor $s\in[1,2]$. On the contrary, not every tuple $(\theta,s)$ corresponds to a $\lambda\in\partial\mathbb{B}_2$, but the tuple describes a $1$-dimensional piecewise linear dynamical system. The cut point criterion for $K$ can be applied to such dynamical systems and, surprisingly, when $s=2$ it recovers the structure of the Douady--Hubbard--Thurston lamination \cite{DouadyHubbard1982,DouadyHubbard1983,DouadyHubbard1984,Thurston2009}.

The use of laminations to study a particular class of connected fractals was carried by Sirvent \cite{Sirvent2008,Sirvent2012,Sirvent2014}. In particular, when $N=2$ the invertible IFS $\{g_j:=f_j\circ s_{\sigma(j)}\}_{j=0}^1$ satisfies the common point property \cite{Sirvent2014}, that is that $\cap_{j=0}^{N-1}g_i(A_\lambda)=\{o_\lambda\}$ and that $g_j^{-1}(o_\lambda)=g_i^{-1}(o_\lambda)\in A_\lambda$ for every $i\neq j$, so Sirvent can construct geodesic laminations associated to $A_\lambda$ based on the information of the symbol space. The group of isometries of the lamination is isomorphic to a subgroup of the full group isometries of the fractal \cite{Sirvent2012}.

\begin{quest}
    Are the two laminations constructions related in any way?
\end{quest}
\begin{quest}
    Can either of the lamination construction be extended to $\Phi_N$ for values of $N>2$?
\end{quest}
Probably the Calegari--Walker construction is more suitable for an extension to $\Phi_N$ since it does not require the common point property which is not satisfied in our case of invertible IFS when $N>2$.

To further motivate the study of the Calegari--Walker laminations, they have numerical evidence (see \cite[Section 7.2]{CalegariWalker2026}) that the set of external angles $\theta$ landing on the $1/3$ bulb of the Mandelbrot set $\mathcal{M}(Q_c)$ corresponds to the tuple $(\theta,2)$. In particular, the associated $\lambda\in\partial\mathbb{B}_2$ resembles a Cantor set, suggesting the existence of some sort of dynamics in the parameter plane.

\section*{Acknowledgement(s)}
I wish to thank Rodrigo P\'{e}rez and Christoph Bandt for helpful discussions. I also wish to thank Bob for taking a chance on me during my undergraduate studies and introducing me to the world of fractals and complex dynamics.

\section*{Disclosure statement}
Generative AI was used to organize and review the content of this article. AI was also used to brainstorm the title of this article and to locate and summarize certain references. I take full responsibility for the mathematical statements, proofs, references and any remaining errors.

\bibliographystyle{plainurl}
\bibliography{references}

\end{document}